\documentclass[aop,MSNbibl,dvips]{arximspdf}

\doi{10.1214/13-AOP856} 
\volume{42}
\issue{2}
\pubyear{2014}
\firstpage{760}
\lastpage{793}

\makeatletter

\newcommand{\underset}[2]{\mathop{#2}_{#1}}

\newcommand{\cal}{\mathcal}

\newtheorem{theo}{Theorem}[section]
\newtheorem{lem}[theo]{Lemma}
\newtheorem{prop}[theo]{Proposition}
\newtheorem{cor}[theo]{Corollary}

\newproclaim{notation}[theo]{Notation}
\newproclaim{remark}[theo]{Remark}
\newproclaim{comment}[theo]{Comment}

\newcommand{\E}{\mathbb{E}}
\renewcommand{\P}{\mathbb{P}}
\newcommand{\B}{{\mathcal B}}
\newcommand{\LL}{{\mathbb L}}
\newcommand{\A}{{\cal P}}
\newcommand{\F}{{\cal F}}
\renewcommand{\H}{{\mathcal H}}
\newcommand{\M}{{\cal M}}

\def\E{\mathbb{E}}
\def\Z{\mathbb{Z}}
\def\R{\mathbb{R}}
\def\B{{\cal B}}
\def\F{{\cal F}}

\makeatother

\begin{document}
\begin{frontmatter}

\title{On martingale approximations and the quenched weak invariance principle}
\runtitle{Martingale approximations in Lp}

\begin{aug}
\author[A]{\fnms{Christophe} \snm{Cuny}\corref{}\ead[label=e1]{christophe.cuny@ecp.fr}}
\and
\author[B]{\fnms{Florence} \snm{Merlev\`{e}de}\ead[label=e2]{florence.merlevede@univ-mlv.fr}}
\runauthor{C. Cuny and F. Merlev\`{e}de}
\affiliation{Ecole Centrale de Paris and Universit\'{e} Paris Est}
\address[A]{Laboratoire MAS\\
Ecole Centrale de Paris\\
Grande Voie des Vignes\\
92295 Chatenay-Malabry cedex\\
France\\
\printead{e1}} 
\address[B]{LAMA (UMR 8050)\\
UPEMLV, CNRS, UPEC\\
Universit\'{e} Paris Est\\
B\^{a}timent Copernic\\
5 Boulevard Descartes\\
77435 Champs-Sur-Marne\\
France\\
\printead{e2}}
\end{aug}

\received{\smonth{2} \syear{2012}}
\revised{\smonth{3} \syear{2013}}

%
\begin{abstract}
In this paper, we obtain sufficient conditions in terms of projective
criteria under which the partial sums of a stationary process with
values in ${\mathcal H}$ (a~real and separable Hilbert space) admits an
approximation, in $\LL^p ({\mathcal H}) $, $p>1$, by a martingale with
stationary
differences, and we then estimate the error of approximation in $\LL^p
({\mathcal H}) $. The results are exploited to further investigate the behavior
of the partial sums. In particular we obtain new projective conditions
concerning the Marcinkiewicz--Zygmund theorem, the moderate deviations
principle and the rates in the central limit theorem in terms of
Wasserstein distances. The conditions are well suited for a large
variety of examples, including linear processes or various kinds of
weak dependent or mixing processes. In addition, our approach suits
well to investigate the quenched central limit theorem and its
invariance principle via martingale approximation, and allows us to
show that they hold under the so-called Maxwell--Woodroofe condition
that is known to be optimal.
\end{abstract}

%

\end{frontmatter}

\section{Introduction}\label{intro}

Since the seminal paper of Gordin \cite{Gor} in 1969, approximation
via a martingale is known to be a nice method to derive limit theorems
for stochastic processes. For instance, the martingale method has been
used successfully by Heyde \cite{He} and Gordin and Lifsic \cite{GL}
to derive central limit \mbox{theorems} for the partial sums of a stationary
sequence, and it has undergone substantial improvements. For recent
contributions where the central limit theory and weak convergence
problems are handled with the help of martingale approximations, let us
mention the recent papers by Maxwell and Woodroofe \cite{MW}, Wu and
Woodroofe~\cite{WW}, Peligrad and Utev \cite{PU}, Merlev\`ede and
Peligrad \cite{MP06}, Zhao and Woodroofe \cite{ZW} and Gordin and
Peligrad \cite{GP}. In all these papers, conditions are then imposed
to be able to implement the martingale method,\vadjust{\goodbreak} namely, to approximate
in a suitable way the partial sums of a stationary process by a
martingale. However, to derive many other kinds of limit theorems from
the martingale method, more precise estimates of the approximation
error of partial sums by a martingale may be useful. We refer to the
recent papers by Wu \cite{Wu}, Zhao and Woodroofe \cite{ZW0}, Cuny
\cite{Cuny1}, Dedecker, Doukhan and Merlev\`ede \cite{DDM} and
Merlev\`ede, Peligrad and Peligrad \cite{MPP} where almost sure
behaviors of the partial sums process have been addressed with the help
of estimates of this approximation error.

In order to say more about these papers and to present our results, let
us first introduce the following notation, giving a way to define
stationary processes.

%
%
\begin{notation}
\label{notation1} Let $(\Omega,{\mathcal{A}},{\mathbb{P}})$ be a probability
space, and let $\theta\dvtx\Omega\mapsto\Omega$ be a bijective
bi-measurable transformation preserving the probability ${\mathbb
{P}}$. Let $\mathcal{F}%
_{0}$ be a $\sigma$-algebra of $\mathcal{A}$ satisfying $\mathcal
{F}_{0}\subseteq\theta^{-1}(\mathcal{F}_{0})$. We then define a
nondecreasing
filtration $(\mathcal{F}_{i})_{i\in\mathbb{Z}}$ by $\mathcal{F}_{i}%
=\theta^{-i}(\mathcal{F}_{0})$, and
a stationary sequence $(X_i)_{i \in\mathbb Z}$ by
$X_i = X_0 \circ\theta^i$ where $X_0$ is a real-valued centered
random variable (or possibly taking values in some real and separable
Hilbert space). The sequence will be called adapted to the filtration
$(\mathcal{F}_{i})_{i\in\mathbb{Z}}$ if $X_0$ is $\mathcal
{F}_{0}$-measurable. Define then the partial sum by $S_n= X_1+X_2+
\cdots+ X_n$. The
following notations will also be used: $\mathcal{F}_{-\infty}=
\bigcap_{i\in\mathbb{Z}}\mathcal{F}_{i}$, $\mathcal{F}_{\infty}=
\bigvee_{i\in\mathbb{Z}}\mathcal{F}_{i}$, $\E_{k}(X)=\E
(X|\mathcal{F}_{k})$, $\A_k(X)= \E_{k}(X)-\E_{k-1}(X)$, and when
$X$ is real-valued, its $\LL^{p} $ norm is denoted by
$\| X \|_{p}= ( \E( | X|^p) )^{1/p}$. We shall also
use the notation $a_{n}\ll b_{n}$ to mean that there exists a numerical
constant $C$
not depending on $n$ such that $a_{n}\leq C b_{n}$, for all positive
integers $n$.
\end{notation}

In all of what follows the sequence $(X_{i})_{i\in\mathbb{Z}}$ is
assumed to be stationary and adapted to $(\mathcal{F}_{i})_{i\in
\mathbb{Z}}$
and the variables are in $\LL^{p}$, for some $p>1$.

In \cite{Wu} and \cite{DDM}, it is assumed that $D = \sum_{i \geq0}
{\mathcal P}_0(X_i)$ converges in $\LL^p$, $p>1$, and estimates of $
\| S_n - M_n \|_p$ where $M_n = \sum_{i=1}^n D \circ\theta^i$
are provided involving either the terms $ \sum_{k \geq n} \|
{\mathcal P}_0(X_k) \|_p$ (see \cite{Wu}) or the terms $ \|\E
_0(S_n) \|_p$ and $\|{\sum_{k \geq n} {\mathcal P}_0(X_k) }\|
_p$; see \cite{DDM}. Those estimates are then exploited to derive
explicit rates in the almost sure invariance principle under projective
conditions that are well adapted to a large variety of examples. The
paper by Merlev\`ede et al. \cite{MPP} addresses different questions
about the almost sure behavior of $S_n$ such as quenched invariance
principles or almost sure central limit theorems. Their proof is based
on a precise estimate of the $\LL^{2} $ approximation error between
the partial sums process and their constructed approximating stationary
martingale, provided that the Maxwell--Woodroofe condition (\ref{MW})
holds. More precisely, in the case where $p=2$, they proved that if
%
%
\begin{equation}\label{MW}
\sum_{k=1}^{\infty}\frac{\|\mathbb{E}_{0}(S_{k})\|_2
}{k^{3/2}}<\infty,
\end{equation}
then there is a martingale $M_n$ with stationary
and square integrable differences such that
%
%
\begin{equation}\label{ApproxMPP}
\| S_{n}-M_{n}\|_2 \ll n^{1/2}
\sum_{k\geq n}\frac{\|
\mathbb{E}_{0}(S_{k})\|_2}{k^{3/2}}. %
\end{equation}
To implement a martingale method for other questions related to the
behavior of the partial sums, as, for instance, rates in the strong laws
of large numbers or in the central limit theorem in terms of
Wasserstein distances, or also moderate deviations principles, the
first question that our paper addresses is the construction of a
stationary martingale $M_n$ in $\LL^{p} $ ($p>1$) in such a way that
an estimate of $\| S_{n}-M_{n}\|_p$ can be given in the spirit
of (\ref{ApproxMPP}). Our Theorem \ref{theoMW} is in this direction.
When $p\geq2$, it states in particular that if
%
%
\begin{equation}\label{MW*p}
\sum_{k=1}^{\infty}\frac{\|\mathbb{E}_{0}(S_{k})\|_p
}{k^{1+1/p}}<\infty,
\end{equation}
then we can construct a stationary sequence $(D_{k}=D \circ\theta
^{k})_{k\in\mathbb{Z}}$ of martingale differences in $\LL^p$ adapted
to $({\mathcal F}_{k})_{k\in\mathbb{Z}}$ such that setting $M_n=\sum
_{k=1}^n D \circ\theta^{k}$,
%
%
\begin{equation}\label{resMW*p}
\| S_{n}-M_{n}\|_p \ll n^{1/2}
\sum_{k\geq[n^{p/2}]}\frac
{\|\mathbb{E}_{0}(S_{k})\|_p
}{k^{1+1/p}}.
\end{equation}
While (\ref{resMW*p}) and (\ref{ApproxMPP}) coincide
when $p=2$, our method of proof is different from the one used in
\cite{MPP}.
In Theorem \ref{theoMW}, we shall consider also the case when $p \in
\ ]1,2[$. The main tools to prove the martingale approximation with the
bound (\ref{resMW*p}) being algebraic computations and Burkholder's
inequality, the estimate also holds for variables taking values in a
separable real Hilbert space. Hence Theorem \ref{theoMW} is stated in
this setting. As we shall see, this martingale approximation result
leads to new projective conditions allowing results concerning the
moderate deviations principle or also estimates of Wasserstein
distances in the CLT; see Sections \ref{submdp} and \ref{subwasser}.
Notice that the projective conditions assumed throughout the paper are
general enough to contain a wide class of dependent sequences.

Another interesting point of our approach and of the approximating
martingale we consider here, is that they lead not only to a useful
estimate of $\| S_{n}-M_{n}\|_p$, but, together with a new
ergodic theorem with rate (see Theorem \ref{theoergo}), they allow
also to show that, under the Maxwell--Woodroofe condition (\ref{MW}),
$\E_0[(S_n-M_n)^2] =o(n)$ $\P$-a.s.; see our Proposition \ref
{propquenched1}. This allows us to give a definitive positive answer to
the question of whether the quenched central limit theorem for
$n^{-1/2}S_n$ holds true under (\ref{MW}). As we shall see, we can
even say more since, using a maximal inequality from Merlev\`ede and
Peligrad \cite{MP}, we establish in Theorem \ref{quenched} that the
functional form of the quenched central limit theorem also holds under
the Maxwell--Woodroofe condition.

Our paper is structured as follows. Section \ref{sectmain} contains
our main results. More precisely, in Section \ref{submartappro} we
construct an approximating martingale with stationary differences in
$\LL^p$ that leads to estimates of the ${\mathbb L}^p$ approximating
error between the partial sums and the constructed martingale; see
Theorem \ref{theoMW}. In Section~\ref{subquenched}, we address the
question of the quenched weak invariance principle under the
Maxwell--Woodroofe condition (\ref{MW}). Section \ref{sectappli} is
devoted to some applications of the estimates given in Theorem \ref
{theoMW} to various kind of limit behavior of the partial sums. In
Section \ref{sectproof}, we prove the results stated in Sections \ref
{submartappro} and \ref{subquenched} and state a new ergodic theorem
with rate (see Theorem \ref{theoergo}) whose proof is postponed in
Appendix \ref{sectergodic}. Some technical results are given and proven
in Appendix \ref{secttechres}.

\section{Main results} \label{sectmain}In complement to Notation
\ref{notation1}, we introduce additional notations used throughout the paper.
%
%
\begin{notation}
\label{notation2} Let ${\mathcal H}$ be a real and separable Hilbert
space equipped with the norm $| \cdot|_{\mathcal H}$. For a random
variable $X$ with values in ${\mathcal H}$, we denote its norm in $\LL
^{p} ({\mathcal H})$ by
$\| X \|_{p,{\mathcal H}}= ( \E( | X|_{{\mathcal H}})^p
)^{1/p}$, and we simply denote $\LL^p=\LL^p(\R)$.
\end{notation}

%
%
\begin{notation}
\label{notation3} Let $p'=\min(2,p)$, $p''=\max(2,p)$ and $q=p''/p'$.
\end{notation}
%
\subsection{Martingale approximation in $\LL^p({\mathcal H})$}
\label{submartappro}
Let $p >1$. In this section, we shall establish conditions in order for
$S_n$ to be approximated by a martingale $M_n$ with stationary
differences in $\LL^p({\mathcal H})$ in such a way that the
approximation error $\| S_n - M_n \|_{p,{\mathcal H}}$ is explicitly
controlled.

Let $(X_n)_{n\in\mathbb{Z}}$ be an adapted stationary sequence in
$\LL^p({\mathcal H})$ in the sense of Notation \ref{notation1}.
When
%
%
\begin{equation}\label{defdiffmart}
D=\sum_{n\geq0} \sum_{k\geq n}
\frac{{\mathcal P}_0 (X_k)}{k+1}
\end{equation}
converges in $\LL^p({\mathcal H})$, then $(D_{k}=D \circ\theta
^{k})_{k\in\mathbb{Z}}$ forms a stationary sequence of martingale
differences in $\LL^p({\mathcal H})$ adapted to $({\mathcal
F}_{k})_{k\in\mathbb{Z}}$. Notice that, by Lemma \ref{lem}, the
series $\sum_{k\geq0} \frac{{\mathcal P}_0 (X_k)}{k+1}$ converges in
$\LL^p({\mathcal H})$ as soon\vspace*{1pt} as $X_0 \in\LL^p({\mathcal H})$. In
addition, note that the series in (\ref{defdiffmart}) converges
in $\LL^p({\mathcal H})$ as soon as the series $ \sum_{k\geq0}
{\mathcal P}_0
(X_k) $ does; see Lemma \ref{lmaabel}.

%
%
\begin{theo}\label{theoMW}
Let $p>1$, and let $(X_n)_{n\in\mathbb{Z}}$ be an adapted stationary
sequence in $\LL^p({\mathcal H})$ in the sense of Notation \ref
{notation1}. Assume that
%
%
\begin{equation}
\label{MWpalpha} \sum_{n\ge1}\frac{\|\E_0(S_n)\|_{p, {\mathcal
H}}}{n^{1+1/p''}}<\infty.
\end{equation}
Then $\sum_{n\ge1}|\sum_{k\ge n} k^{-1} {\mathcal
P}_0(X_{k-1})|_{{\mathcal H}}$ converges in
$\LL^p$ and setting $M_n=\break \sum_{k=1}^{n}D\circ\theta^k$ where $D$ is
defined by (\ref{defdiffmart}), the following inequality holds:
%
%
\begin{equation}
\label{rate} \|S_n-M_n\|_{p, {\mathcal H}} \ll
n^{ 1/p' } \sum_{k\ge[n^{q}]}\frac{\|\E_0(S_k)\|_{p, {\mathcal
H}}}{k^{1+1/p''}}.
\end{equation}
%
\end{theo}
%
%
\begin{remark}\label{remMP}
Let $p>1$ and $\alpha\in\ ]0, 1/p'']$. Let us introduce the following
assumption:
%
%
\begin{equation}
\label{MWpalpha*} \sum_{n\ge1}\frac{\|\E_0(S_n)\|_{p,{\mathcal
H}}}{n^{1+\alpha
}}<\infty.
\end{equation}
Assume that (\ref{MWpalpha*}) holds with $\alpha=\min( 1/2, 2/p^2)$.
By combining (\ref{rate}) with Corollary~22 of \cite{MP} (with the
norm $|\cdot|_{{\mathcal H}}$ replacing the absolute values) we
have
%
%
\begin{equation}
\label{maxbis} \Bigl\|\max_{1\le k\le n}|S_k-M_k|_{{\mathcal H}}
\Bigr\|_{p} =o\bigl(n^{1/p}\bigr).
\end{equation}
Notice also that if $p >2$ and (\ref{MWpalpha*}) holds with $\alpha
\in\ ]2/p^2, 1/p]$, then (\ref{rate}) combined with the maximal
inequality (7) of \cite{MP} (with the norm $|\cdot|_{{\mathcal H}}$
replacing the absolute values) implies that
\[
\Bigl\|\max_{1\le k\le n}|S_k-M_k|_{{\mathcal H}}
\Bigr\|_p =o\bigl(n^{\alpha p/2}\bigr).
\]
The fact that the maximal inequality (7) of \cite{MP} is still valid
when the variables take values in a Hilbert space comes from the fact
that its proof is only based on chaining arguments (still valid in
functional spaces by replacing the absolute values by the corresponding
norms) and on Doob's maximal inequality that also holds in Hilbert
spaces. Since Corollary 22 of \cite{MP} is proved via their maximal
inequality (7), it is still valid in the Hilbert space setting.
\end{remark}
\subsection{Martingale approximation under $\mathbb{P}_0$ and the quenched
(weak) invariance principle} \label{subquenched}

Limit theorems for stochastic processes that do not start from
equilibrium are timely and motivated by evolutions in a quenched random
environment. Recent discoveries by Voln\'{y} and Woodroofe \cite{VW}
show that many of the central limit theorems satisfied by classes of
stochastic processes in equilibrium fail to hold when the processes are
started from a point. In this section, we address the question of
whether the Maxwell--Woodroofe condition (\ref{MW}) is sufficient for
the validity of the quenched central limit theorem since this condition
is known to be optimal; see, for example, \cite{PU} or
\cite{Volny} where the optimality of this condition is discussed. This
question starts with a result in Borodin and Ibragimov (\cite{BI},
Chapter 4) stating that
if $\|\mathbb{E}_{0}(S_{n})\|_2$ is bounded, then one has the
CLT starting at a point in its
functional form. Later, works by Derriennic and Lin
(see \cite{DL,DL2,DL3}), Zhao
and Woodroofe \cite{ZW0}, Cuny and Lin \cite{CL}, Cuny \cite{Cuny1}
and Merlev\`ede, Peligrad and Peligrad \cite{MPP} improved on this
result by imposing weaker and weaker conditions on $\|\mathbb{E}%
_{0}(S_{n})\|_2$, but always stronger than (\ref{MW}). Let us
mention that a result in Cuny and Peligrad \cite{CP} shows
that the condition $\sum_{k=1}^{\infty}\|\mathbb
{E}_{0}(X_{k})\|_2
/k^{1/2}<\infty$ is sufficient for the quenched CLT. It is also
sufficient for the quenched weak invariance principle
by a recent result of Cuny and Volny \cite{CV}.

As we shall see in the proof of Theorem \ref{quenched} below, the
approximating martingale that we defined in Section \ref{submartappro}
also allows us to show that, under (\ref{MW}), $\lim_{n \rightarrow
\infty} n^{-1} \E_0(|S_n-\E_0 (S_n) - M_n|^2) = 0 $ $\P$-a.s.
Combined with a new ergodic theorem with rate (see our Theorem \ref
{theoergo}) and a maximal inequality from Merlev\`ede and Peligrad
\cite{MP}, this implies that the quenched CLT in its functional form
holds under the Maxwell--Woodroofe condition (\ref{MW}).

To state that result we need some further notations.
Let us first assume the existence of a regular version of the conditional
probability on ${\mathcal A}$ given $\F_0$; that is, we assume the existence
of a transition probability $K(\cdot,\cdot)$ on $(\Omega,{\mathcal
A})$, such that
for every $A \in{\mathcal A}$, $K(\cdot, A)$ is a version of $\E
({\mathbf1}_A|\F_0)$.
Then we denote by $\E_\omega$ the expectation with respect
to $K(\omega,\cdot)$. We also define the Donsker
process $W_n$ by $W_n(t)=n^{-1/2}(S_{[nt]}+(nt-[nt])X_{[nt]+1})$.

%
%
\begin{theo}\label{quenched}
Let $(X_n)_{n\in\mathbb{Z}}$ be an adapted stationary sequence\break  in
$\LL^2$ in the sense of Notation \ref{notation1}. Assume
that (\ref{MW}) holds.
Then\break  $\sum_{n\ge1}|\sum_{k\ge n} k^{-1} {\mathcal P}_0(X_{k-1})|$
converges in
$\LL^2$, and setting $M_n=\sum_{k=1}^{n}D\circ\theta^k$ where $D$
is defined by (\ref{defdiffmart}), the following holds:
%
%
\begin{equation}
\label{resquen} \frac{\E_0(\max_{1\le k \le n}|S_k-M_k|^2)
}{n}\underset{n\rightarrow+\infty} {
\longrightarrow} 0 \qquad\mbox{$\P$-a.s.}
\end{equation}
In particular, $(S_n)$ satisfies the following quenched weak invariance
principle: there exists $\Omega_0\in{\mathcal A}$ with $\P(\Omega
_0)=1$ such that for every $\omega\in\Omega_0$, for any continuous
and bounded function $f$ from $%
(C([0,1]), \|\cdot\|_{\infty})$ to ${\mathbb{R}}$,
%
%
\begin{equation}
\label{FCLT} \lim_{n \rightarrow\infty}\mathbb{E}_{\omega}
\bigl(f(W_n)\bigr)=\int f \bigl(z \sqrt{\eta(\omega) }\bigr)W(dz),
\end{equation}
where $\eta= \lim_{n \rightarrow\infty} n^{-1}\E(S_n^2 | {\mathcal I})
= \lim_{n \rightarrow\infty} n^{-1}\E_0 (S_n^2 )$ in $\LL^1$, and
$W$ is the distribution of a standard Wiener process. Here ${\mathcal
I}$ is the invariant sigma field, that is, ${\mathcal I}=\{A\in
{\mathcal A}\dvtx \theta^{-1}(A)=A\}$.
\end{theo}
It follows from Comment \ref{commentrho} that if the $\rho$-mixing
coefficients of $(X_n)_{n\in\mathbb{Z}}$ satisfy $\sum_{k \geq0}
\rho(2^k) < \infty$, then the quenched invariance principle holds.
Hence the CLT from Ibragimov \cite{Ib} for $\rho$-mixing sequences
that is known to be essentially optimal, is also quenched.

A careful analysis of the proof of Theorem \ref{quenched}
shows that if the
random variables are assumed to be in $\LL^2(\H)$, then under (\ref
{MWpalpha}) with $p=2$,
the almost sure convergence (\ref{resquen}) still holds with the norm
$| \cdot|_{\H}$
replacing the absolute values.

Theorem \ref{quenched} has an interesting interpretation in the
terminology of additive
functionals of Markov chains.
Let $(\xi_{n})_{n\geq0}$ be a Markov chain with values in a Polish
space $S$,
so that there exists a regular transition probability $P_{\xi_1|\xi
_0=x}$. Let $P$
be the transition kernel defined by $%
P(g)(x )=P_{\xi_1|\xi_0=x}(g)$ for any bounded measurable function
$g$ from $S$ to ${\mathbb R}$, and assume that
there exists an invariant probability $\pi$ for this transition kernel,
that is, a probability measure on $S$ such that $\pi(g)=\pi(P(g))$
for any
bounded measurable function $g$ from $S$ to ${\mathbb R}$.
Let then $\mathbb{L}^{2}(\pi)$ be
the set of functions from $S$ to ${\mathbb R}$ such that $\pi(g^2)
<\infty$.
For $g\in\mathbb{L}^{2}(\pi)$ such that $\pi(g)=0$, define $%
X_{i}=g(\xi_{i})$. In this setting condition (\ref{MW}) is $\sum_{n\geq
1} n^{-3/2} \|\sum_{k=1}^n P^k (g) \|_{\mathbb
{L}^{2}(\pi)}
<\infty$. In the context of a Markov chain, the conclusion of Theorem
\ref{quenched} is also known under
the terminology of functional CLT started at a point. To rephrase it,
let $\mathbb{P}%
^{x}$ be the probability associated to the Markov chain started from
$x$ and
let $\mathbb{E}^{x}$ be the corresponding expectation. Then, for $\pi
$-almost every $x \in S$, for any continuous and bounded function $f$
from $(C([0,1]), \|\cdot\|_{\infty})$ to~${\mathbb{R}}$,\looseness=1
%
%
\begin{equation}
\label{quenchedMC} \lim_{n \rightarrow\infty} \mathbb{E}^x
\bigl(f(W_n)\bigr) = \int f (z \sqrt{\eta_x})W(dz),
\end{equation}\looseness=0
where $\eta_x:=\lim_n \E^x(S_n^2)/n$. Note that Theorem \ref
{quenched} improves Corollary 5.10 of \cite{Cuny1} stated for Markov
chains with normal Markov operator. Let us mention that convergence
(\ref{quenchedMC}) has also been obtained recently in Dedecker,
Merlev\`ede and Peligrad \cite{DMP} under the condition $\sum_{k\geq
0} \pi(|gP^{k}(g)|) <\infty$. The latter condition and (\ref{MW})
are of independent interests; see Section 5.2 of~\cite{DMP}.

\section{Applications} \label{sectappli}
As we mentioned in the \hyperref[intro]{Introduction}, having estimates
of the
approximation error of partial sums by a martingale can be useful to
derive different kinds of limit theorems for the partial sums
associated with a stationary process. For instance, starting from
(\ref{ApproxMPP}), Merlev\`ede et al. \cite{MPP} have obtained
sufficient projective conditions in order for the partial sums to
satisfy either the law of the iterated logarithm or the almost sure
central limit theorem. In this section, we shall use our estimate
(\ref{rate}), either to give new projective conditions under which the
partial sums associated with a stationary process satisfy a moderate
deviations type results, or to analyze the rates of convergence in the
CLT in terms of Wasserstein distances. Before
stating those results we provide a simple and direct application of our results,
leading to new projective criteria to obtain rates in the SLLN.

\subsection{Strong law of large numbers with rate} \label{subslln}

Our martingale approximation in $\LL^p$ for $1<p<2$ combined with our
new ergodic theorem with rate (see Theorem~\ref{theoergo}) allows us to
derive very directly a projective condition for the
Marcinkiewicz--Zygmund strong law of large numbers.

%
%
\begin{theo}\label{MZ}
Let $1<p<2$, and let $(X_n)_{n\in\mathbf Z}$ be an adapted stationary
sequence in $\LL^p({\mathcal H})$ in the sense of Notation \ref{notation1}.
Assume that
\[
\sum_{n\ge2} \log n
\frac{\|\E_0(S_n)\|_{p,{\mathcal H}}}{n^{3/2}}<\infty.
\]
Then there exists a stationary martingale $(M_n)_{n \geq1}$ in $\LL
^p({\mathcal H})$, such that \mbox{$|S_n -M_n|_{{\mathcal H}} = o(n^{1/p})$}
${\mathbb P}$-a.s. In particular, we have $\vert S_n \vert_{\H
}=o(n^{1/p})$ ${\mathbb P}$-a.s.
\end{theo}

\begin{pf}
Using Theorem \ref{theoergo}, the first part of the result will follow
if we can prove that
$\sum_{n\ge1}n^{-1-1/p}\| S_n -M_n\|_{p,{\mathcal H}} <\infty$. This
convergence follows by using Theorem \ref{theoMW} to control $\|S_n
-M_n\|_{p,{\mathcal H}}$. For the last part of the theorem, it suffices
to notice that by the Marcinkiewicz--Zygmund strong law of large
numbers for martingales $|M_n|_{{\mathcal H}}=o(n^{1/p})$
${\mathbb P}$-a.s. for any $p \in\ ]1,2[$ as soon as the martingales
are in $\LL
^p({\mathcal H})$; see Woyczy\'nski \cite{WO}.
\end{pf}

\subsection{Moderate deviations} \label{submdp}

The aim of this section is to obtain asymptotic expansions for
probabilities of moderate deviation for stationary adapted real-valued
processes under projective criteria; more precisely we want to study
the asymptotic behavior of $\P( S_n \geq\sigma\sqrt n r_n)$ where
$(r_n)$ is a sequence of positive numbers that diverges to infinity at
an appropriate rate and $\sigma=\lim_{n \rightarrow\infty}\|S_n \|
_2/\sqrt n$. Specifically, we aim to
find the zone for $x$ of the following moderate deviations principle:
%
%
\begin{equation}\label{moderate}
{\frac{{\mathbb{P}(S_{n}\geq x\sigma\sqrt n r_n)}}{1-\Phi(x
r_n)}}%
=1+o(1),%
\end{equation}
where $\Phi(x)$ is the standard normal distribution function. If
$r_n=r>0$ is
fixed, then (\ref{moderate}) is essentially the well-known central
limit theorem.
However, for the case when $r=r_{n}$ is allowed to tend to infinity, the
problem of \textit{moderate deviation} probabilities is to find all
the possible speed of convergence of $r_{n}\rightarrow\infty$ such that
(\ref{moderate}) holds. It is a challenging problem to establish
moderate deviations principle (MDP) for dependent variables. However,
starting from the deep results of Grama \cite{GR} and of Grama and
Haeusler \cite{GH} for martingales, Wu and Zhao \cite{WZ} showed that
it is possible to obtain MDP results for a certain class of stationary
processes such as functions of an i.i.d. sequence as soon as the
partial sum process can be well approximated by a martingale. Using our
Theorem \ref{theoMW}, we shall give sufficient conditions for the MDP
to hold that are
different than those obtained by Wu and Zhao~\cite{WZ}.

Let us first start with some notation and definitions.

Let $p \in(2,4]$. For $x >1$, let $r_x >0$ be the solution to the equation
\[
x=(1 + r_x)^{\nu(p)} \exp\bigl(r_x^2/2
\bigr) \qquad\mbox{where } \nu(p) = \cases{ p+1, &\quad if $2< p \leq3$,
\cr
3p-3, &\quad if $3 <
p \leq4$.}
\]
The function $\nu(p)$ results from the martingale MDP as obtained in
\cite{GR} and in \cite{GH}; see also Theorem 2 and Remark 5 in \cite
{WZ}. In addition, by Remark 1 in \cite{GH}, as $x \rightarrow\infty
$, $r_x$ has the asymptotic expansion $r_x^2=2 \log x - 2 [ \nu( p)
+o(1)] \log( 1 + \sqrt{2 \log x} )$.

Let $\tau_n \rightarrow\infty$ be a positive sequence of numbers and
$(U_n)$ a sequence of real valued random variables such that $U_n
\rightarrow^{{\mathcal D}}
{\mathcal N} (0,1)$. We shall say that $(U_n)$ satisfies the moderate
deviation principle (MDP) with rate $\tau_n$ and exponent $p>0$ if for
every $a>0$ there exists a positive constant $C=C_{a,p}$ depending
neither on $x$ nor on $n$ such that
\[
\max\biggl\{ \biggl| \frac{\P(U_n \geq r_x)}{ 1 - \Phi(r_x)} -1 \biggr|, \biggl| \frac{\P
(U_n \leq-r_x)}{ 1 - \Phi(-r_x)} -1 \biggr| \biggr\}
\leq C \biggl( \frac{x}{\tau_n} \biggr)^{1/(1+p)}
\]
holds uniformly in $x \in[1, a \tau_n]$. Therefore $\tau_n$ gives a
range for which the MDP holds.

%
%
\begin{theo}\label{theoMDP}
Let $2 <p \le4$, and let $(X_n)_{n\in\mathbb{Z}}$ be an adapted
stationary sequence in $\LL^p$ in the sense of Notation \ref
{notation1}. Assume that
%
%
\begin{equation}
\label{condMDP1}\sum_{n\ge1}\frac{\|\E_0(S_n)\|
_p}{n^{1+2/p^2}}<\infty
\quad\mbox{and}\quad\sum_{n\ge1} \frac{1}{n^{2/p}} \sum
_{k \geq n} \frac{\|\E_{-n}(S_k)\|_2}{k^{3/2}} <\infty.
\end{equation}
Assume in addition that
%
%
\begin{equation}
\label{condSnMDPnondya}\sum_{n\ge1}
\frac{1}{n^{1+2/p}}\bigl\|\E_{-n}\bigl(S_n^2\bigr) -
\E\bigl(S_n^2\bigr)\bigr\|_{p/2}<\infty.
\end{equation}
Then $n^{-1}\E(S_n^2)$ converges to some nonnegative number $\sigma
^2$, and if $\sigma>0$, $ ( \frac{S_n}{\sigma\sqrt n} )_{n
\geq1} $ satisfies the MDP with rate $\tau_n=n^{p/2-1}$ and exponent $p$.
\end{theo}
\begin{pf}
Analyzing the proof of Theorem 1 in \cite{WZ}, we infer
that the theorem will be proven if we can show that there exists a $\LL
^p$ stationary sequence $(D_i)_{i \in{\mathbb Z}}$ of martingale
differences with respect to $({\mathcal F}_i)_{i \in{\mathbb Z}}$ such
that setting $M_n=\sum_{i=1}^n D_i$,
%
%
\begin{equation}
\label{condRnMDP}\| S_n -M_n \|_p= o
\bigl(n^{1/p}\bigr)
\end{equation}
and
%
%
\begin{equation}
\label{condMnMDP1} \Biggl\|\sum_{i=1}^n
\E_{i-1}\bigl(D_i^2\bigr) -\E
\bigl(D_i^2\bigr) \Biggr\|_{p/2}= O
\bigl(n^{2/p} \bigr).
\end{equation}
According to Theorem \ref{theoMW} combined with Remark \ref{remMP},
the first part of condition (\ref{condMDP1}) implies (\ref
{condRnMDP}). On the other hand, since $1<p/2 \leq2$, according to
Theorem 3 in \cite{WZ} applied to the stationary sequence $(\E_{i-1}(D_i^2)
-\E(D_i^2))_{i \geq1}$ and using the fact that $M_n$ is a martingale,
(\ref{condMnMDP1}) holds if
%
%
\begin{equation}
\label{condMnMDPdyadicwithoutgap}\sum_{k\ge0}
\frac{1}{2^{2k/p}}\bigl\| \E_0\bigl(M_{2^k}^2\bigr) -
\E\bigl(M_{2^k}^2\bigr)\bigr\|_{p/2}<\infty.
\end{equation}
We notice now that since $M_n$ is a stationary martingale, for any $r
\geq1$,
%
%
\begin{eqnarray}
\label{compdyanorm}
&&
\bigl\|\E_0\bigl(M_{2^k}^2
\bigr) -\E\bigl(M_{2^k}^2\bigr) \bigr\|_{r}\nonumber\\
&&\qquad =
\Biggl\|\sum_{i=0}^{k-1} \bigl( \E
_{-2^i}\bigl(M_{2^i}^2\bigr) -\E
\bigl(M_{2^i}^2\bigr) \bigr) \circ\theta^{2^i}+
\bigl( \E_0\bigl(D_1^2\bigr) -\E
\bigl(D_1^2\bigr) \bigr) \Biggr\|_{r}
\nonumber\\[-8pt]\\[-8pt]
&&\qquad \leq \sum_{i=0}^{k-1} \bigl\|
\E_{-2^i}\bigl(M_{2^i}^2\bigr) -\E
\bigl(M_{2^i}^2\bigr) \bigr\|_{r} + \bigl\|
\E_0\bigl(D_1^2\bigr) -\E\bigl(D_1^2
\bigr) \bigr\|_{r}
\nonumber\\
&&\qquad \leq 2 \sum_{i=1}^{k-1} \bigl\|
\E_{-2^i}\bigl(M_{2^{i-1}}^2\bigr) -\E
\bigl(M_{2^{i-1}}^2\bigr) \bigr\|_{r} + 2 \bigl\|
\E_0\bigl(D_1^2\bigr) -\E\bigl(D_1^2
\bigr) \bigr\|_{r}.\nonumber
\end{eqnarray}
It follows that (\ref{condMnMDPdyadicwithoutgap}) is equivalent to
$\sum_{k\ge0} 2^{-2k/p}\|\E_{-2^{k+1}}(M_{2^k}^2)
-\E(M_{2^k}^2)\|_{p/2}<\infty$. Due to the subadditivity of the
sequence $(\|\E_{-2n}(M_{n}^2)
-\E(M_{n}^2)\|_{p/2})_{n \geq1}$, the latter condition is equivalent to
%
%
\begin{equation}
\label{condMnMDP}\sum_{n\ge1} \frac{1}{n^{1+2/p}} \bigl\|
\E_{-2n}\bigl(M_n^2\bigr) -\E
\bigl(M_n^2\bigr)\bigr\|_{p/2}<\infty;
\end{equation}
see Lemma 2.7 in \cite{PU}. Using now Proposition \ref{passMnSn}, we
infer that (\ref{condMnMDP}) holds if (\ref{condSnMDPnondya}) and the
second part of (\ref{condMDP1}) hold and if $\sum_{n\ge
1}n^{-(1+4/p^2)} \|\E_0(S_n)\|^2_p<\infty$.\vspace*{1pt} To end the proof, it
suffices to notice that since $(\|\E_0(S_n)\|_p)_{ n \geq1}$ is a
subadditive sequence, the latter condition is satisfied provided the
first part of (\ref{condMDP1}) is satisfied as well; see item 3 of
Lemma 37 in \cite{MP}.
\end{pf}

The quantities involved in conditions (\ref{condMDP1}) and (\ref
{condSnMDPnondya}) can be handled by controlling norms of individual
summands which involve terms such as $\E_0(X_iX_j) $ and $ \E_0(X_i)
$. The latter quantities can be then in turn controlled by using
various mixing or dependence coefficients; see, for example, \cite
{DDM}. For instance, as a corollary of Theorem \ref{theoMDP}, the
following result holds; its proof is omitted since it follows the lines
of the proof of Corollary 2.1 in \cite{DDM}.
%
%
\begin{cor}
Let $2 <p \le4$, and let $(X_n)_{n\in\mathbb{Z}}$ be an adapted
stationary sequence in $\LL^p$ in the sense of Notation \ref
{notation1}. Assume that there exists $\gamma\in\ ]0,1]$ such that
\[
\sum_{n>0} \frac{n^{(p-2)/(\gamma p)}}{ n^{1/p}} \bigl\|
\E_0(X_n) \bigr\|_p < \infty
\]
and
\[
\sum
_{n>0} \frac{n^{\gamma}}{ n^{2/p}} \sup_{i \geq j \geq n }
\bigl\| \E_0(X_iX_j) - \E(X_iX_j)
\bigr\|_{p/2} < \infty.
\]
Then the conclusion of Theorem \ref{theoMDP} holds with $\sigma^2 =
\sum_{k \in\Z} \operatorname{Cov}(X_0,X_k)$.
\end{cor}
As in \cite{DDM}, this result may be used, for instance, to derive
under which conditions the partial sum of a function $f$ of the stationary
Markov chain $(\xi_k)_{k \in{\mathbb Z}}$ with transition $Kf (x)
=\frac{1}{2}( f (x+a) + f (x-a) )$, when $a$ is irrational in $[0,1]$
and badly approximable by rationals, satisfy the conclusion of Theorem
\ref{theoMDP}.
For instance, one can prove that if $f$ is three times differentiable,
$ ( \frac{S_n(f)}{\sigma(f) \sqrt n} )_{n \geq1} $ satisfies
the MDP with rate $\tau_n=n$ and exponent $4$ provided that $\sigma
(f)>0$. Here $S_n (f) = \sum_{k=1}^n (f(\xi_k) - m(f) )$
where $m$ is the Lebesgue--Haar measure and $\sigma^2(f) =
m((f-m(f))^2) + 2 \sum_{n>0} m( f K^n(f-m(f)))$.

Since in Theorem \ref{theoMDP} the conditions are expressed in terms
of the conditional expectation of the partial sum or of its square, it
is also possible to obtain applications for mixing sequences. As an
example, the following corollary gives conditions in terms of $\rho
$-mixing coefficients as defined in Comment~\ref{commentrho}.

%
%
\begin{cor} \label{corrhoMDP}
Let $2<p\le4$, and let $p\le\alpha\le4$. Let $(X_n)_{n\in\mathbb
{Z}}$ be an adapted stationary sequence in $\LL^{\alpha}$ in the
sense of Notation \ref{notation1}. Let $(\rho(n))_{n \geq1}$ be its
associated rho-mixing coefficients as defined in (\ref{defrho}).
Assume that
%
%
\begin{equation}
\label{condmdprho} \sum_{n \geq1} \frac{\rho^{2/p}(n)}{n^{1/2+2/p^2}} <
\infty\quad\mbox{and}\quad\sum_{n \geq1} \frac{\rho^{s}(n)}{n^{2/p}} <
\infty\qquad\mbox{where $s = 2(\alpha-2)/\alpha$}.\hspace*{-25pt}
\end{equation}
Then the conclusion of Theorem \ref{theoMDP} holds with rate $\tau
_n=n^{p/2-1}$ and exponent $p$.
\end{cor}
Notice that if $\alpha= 4$, condition (\ref{condmdprho}) reduces to
its first part.

\begin{pf*}{Proof of Corollary \ref{corrhoMDP}}
Let us prove that the first part of (\ref
{condMDP1}) holds. With this aim, we first notice that, due to the
subadditivity of the sequence $( \|\E_0 (S_n) \|_{p})_{n \geq
1}$, this condition is equivalent to (see Lemma 2.7 in~\cite{PU})
%
%
\begin{equation}
\label{firstcondMDP1dya} \sum_{k\ge0}
\frac{\|\E_0(S_{2^k})\|_{p}}{2^{2k/p^2}}<\infty.
\end{equation}
Since $p>2$, (\ref{condmdprho}) implies that $\sum_{k \geq0} \rho
^{2/p}(2^k) < \infty$. Therefore, by using (\ref{majrho}), it follows
that (\ref{firstcondMDP1dya}) is satisfied as soon as $\sum_{k\ge
0}2^{-2k/p^2} \sum_{i=0}^k 2^{i/2} \rho^{2/p}(2^i) < \infty$, which
is equivalent to the first part of
condition (\ref{condmdprho}).

We prove now that the second part of (\ref{condMDP1}) holds. Due to
the monotonicity of the sequence $ ( \sum_{\ell\geq n} \ell
^{-3/2} \|\E_{-n}(S_\ell)\|_2 )_{n \geq1}$, the second part of
(\ref{condMDP1}) is equivalent to
%
%
\begin{equation}
\label{condMDP12rewritten} \sum_{k\ge0}
\frac{2^k}{2^{2k/p}} \sum_{j \geq k} 2^{-3j/2} \sum
_{\ell= 2^j}^{2^{j+1}-1} \bigl\|\E_{-2^k}(S_{\ell})
\bigr\|_2 < \infty.
\end{equation}
To prove the above condition, we first notice that by stationarity, for
any $\ell\in\{2^j,\ldots, 2^{j+1}-1 \}$,
\begin{eqnarray*}
\bigl\|\E_{-2^k}(S_{\ell})\bigr\|_2 &\leq& \bigl\|
\E_{-2^k}(S_{\ell} -S_{2^j})\bigr\|_2 + \bigl\|
\E_{-2^k}(S_{2^j})\bigr\|_2
\\
& \leq& \bigl\|\E_{-2^k-2^j}(S_{\ell-2^j} )\bigr\|_2 + \sum
_{s=0}^{j-1}\bigl\|\E_{-2^k-2^s}(S_{2^s})
\bigr\|_2 \\
&&{}+ \bigl\|\E_{-2^k}(X_1)\bigr\|_2.
\end{eqnarray*}
Since, for any positive integers $r$ and $t$, $\| \E_{-r}(S_t)\|_2 \ll
\rho(r)\sqrt t$, it follows that
\[
\sum_{\ell= 2^j}^{2^{j+1}-1} \bigl\|\E_{-2^k}(S_{\ell})
\bigr\|_2 \ll2^{3j/2} \rho\bigl(2^j\bigr) +
2^j \rho\bigl(2^k\bigr) + 2^j \sum
_{s=0}^{j-1} 2^{s/2}\rho\bigl(2^k
+2^s\bigr).
\]
So overall, since $p >2$, we infer that
%
%
\begin{equation}
\label{corrhomdpp1bis} \sum_{k\ge0}\frac{2^k}{2^{2k/p}}
\sum_{j \geq k} 2^{-3j/2} \sum
_{\ell= 2^j}^{2^{j+1}-1} \bigl\|\E_{-2^k}(S_{\ell})
\bigr\|_2 \ll\sum_{k\ge0} 2^{k(1-2/p)}\rho
\bigl(2^k\bigr).
\end{equation}
Noticing that (\ref{condmdprho}) implies in particular that
%
%
\begin{equation}
\label{implicationcondmdprho} \rho\bigl(2^k\bigr) = o
\bigl(2^{-k(p^2-4)/(4p)} \bigr) \qquad\mbox{as }k \rightarrow\infty,
\end{equation}
and taking into account that $p >2$, we then infer that the sums in the
right-hand side of (\ref{corrhomdpp1bis}) are finite under (\ref
{condmdprho}). This ends the proof of (\ref{condMDP12rewritten}), hence
the second part of (\ref{condMDP1})
holds.

It remains to show that (\ref{condSnMDPnondya}) is satisfied. Note
first that since $p \in\ ]2,4]$ and $\alpha\geq p$,
\begin{eqnarray*}
\bigl\|\E_{-n}\bigl(S^2_{n}\bigr) -\E
\bigl(S_{n}^2\bigr)\bigr\|_{p/2} &\leq&\bigl\|
\E_{-n}\bigl(S^2_{n}\bigr) -\E
\bigl(S_{n}^2\bigr)\bigr\|_{ \alpha/2}\\
&\leq&\sup
_{Z \in B^{ \alpha/(\alpha
-2)}({\mathcal F}_{-n})} \operatorname{Cov} \bigl(Z, S^2_{n}\bigr),
\end{eqnarray*}
where ${B}^r({\mathcal F}_{-n})$ stands for the set of ${\mathcal
F}_{-n}$-measurable
random variables such that $\|Z\|_r \leq1$. Using then Theorem 4.12 in
\cite{BR}, we get that
\[
\bigl\|\E_{-n}\bigl(S^2_{n}\bigr) -\E
\bigl(S_{n}^2\bigr)\bigr\|_{p/2} \leq2^{1-s}
\rho^s (n) \bigl\| S^2_{n}\bigr\|_{\alpha/2}
= 2^{1-s} \rho^s (n) \| S_{n}
\|^2_{\alpha},
\]
where $s = 2(\alpha-2)/\alpha$.
Now the first part of (\ref{condmdprho}) implies $\sum_{k >0} \rho
^{1/2} (2^k) < \infty$ [see also (\ref{implicationcondmdprho})],
therefore $ \| S_{n}
\|_{ \alpha} \ll n^{1/2}$; see \cite{Pe} or \cite{Sh}. Hence
%
%
\begin{equation}
\label{corrhomdpp3} \bigl\|\E_{-n}\bigl(S^2_{n}
\bigr) -\E\bigl(S_{n}^2\bigr)\bigr\|_{p/2} \ll n
\rho^s (n),
\end{equation}
which proves that (\ref{condSnMDPnondya}) holds as soon as the second
part of (\ref{condmdprho}) does. This ends the proof of the corollary.
\end{pf*}

\subsection{Rates of convergence for Wasserstein distances in the
CLT} \label{subwasser}

Let ${\mathcal L}(\mu, \nu)$ be the set of probability laws on
$\mathbb R^2$ with marginals $\mu$ and $\nu$. Let us consider the
Wasserstein distances
of order $r\geq1$ defined by
\[
W_r(\mu, \nu) = \inf\biggl\{ \biggl(\int|x-y|^r P(dx,
dy) \biggr)^{1/r}\dvtx P \in{\mathcal L}(\mu, \nu) \biggr\}.
\]
Let $p\in\ ]2,3[$, and let $(X_n)_{n\in\mathbb{Z}}$ be an adapted
stationary sequence in ${\mathbb L}^p$ in the sense of Notation \ref
{notation1}. Denote by $P_{S_{n}/n^{1/2}}$ the law of $S_{n}/n^{1/2}$
and by
$G_{\sigma^2}$ the normal distribution $ {\mathcal N}(0, \sigma^2)$
where $\sigma^2=\lim_{n \rightarrow\infty}n^{-1}\E(S_n^2)$
provided the limit exists. Starting from Theorem 2.1 in \cite{DMRi}
and using our Theorem \ref{theoMW}, we get the following result
concerning the order of $W^r_{r}( P_{n^{-1/2}S_n},G_{\sigma^2}
)$ where $r \in[1,p]$.
%
%
\begin{theo}\label{Wasser}
Let $2<p \le3$ and let $1\le r \le p$. Let $(X_n)_{n\in\mathbb{Z}}$
be an adapted stationary sequence in $\LL^p$ in the sense of Notation
\ref{notation1}. Assume that (\ref{condSnMDPnondya}) holds and that
%
%
\begin{equation}
\label{Cond1carre} \sum_{n \geq1}\frac{1}{n^{3-p/2}}
\bigl\|\E_{-n} \bigl(S_{n}^2 \bigr) - \E
\bigl(S_{n}^2 \bigr)\bigr\|_{1 + \gamma} < \infty\qquad\mbox{for some
$\gamma>0$}.
\end{equation}
Assume in addition that
%
%
\begin{equation}
\label{cond1Wasser} \sum_{n \geq1} \frac{\|\E_0(S_n) \|
^2_{p}}{n^{1+4/p^2}}<
\infty
\end{equation}
and that
%
%
\begin{eqnarray}
\label{cond2Wasser} \sum_{n\ge1}\frac{\|\E_0(S_n)\|_2}{n^{(5-p)/2}}&<&
\infty\qquad\mbox{if $r \in[1,2]$} \quad\mbox{and}\nonumber\\[-8pt]\\[-8pt]
\bigl\|\E_0(S_n)\bigr\|_r&=&O\bigl(n^{(3-p)/r}\bigr) \qquad\mbox{if $r \in\
]2,p]$}.\nonumber
\end{eqnarray}
Then $n^{-1}\E(S_n^2)$ converges to some nonnegative number $\sigma
^2$, and $
W^r_{r}( P_{n^{-1/2}S_n},\break  G_{\sigma^2}
)=O(n^{1-p/2}) $.
\end{theo}
The above result improves Theorem 3.1 in Dedecker, Merlev\`ede and Rio
\cite{DMRi} that imposes the series $\sum_{n>0}\E(X_n | {\mathcal
F}_0) $ to converge in ${\mathbb L}^p$ instead of the weaker conditions
(\ref{cond1Wasser}) and (\ref{cond2Wasser}).

When $\rho$-mixing sequences are considered, applying Theorem \ref
{Wasser} we derive the following corollary (its proof is omitted since
it uses similar bounds as those obtained in the proof of Corollary \ref
{corrhoMDP}).

%
%
\begin{cor} \label{corrhowasser}
Let $2<p\le3$ and let $p\le\alpha\leq4$. Let $(X_n)_{n\in\mathbb
{Z}}$ be a adapted stationary sequence in $\LL^{\alpha}$ in the sense
of Notation \ref{notation1}. Let $(\rho(n))_{n \geq1}$ be its
associated rho-mixing coefficients as defined in (\ref{defrho}).
Assume that
\[
\sum_{n \geq1} \frac{\rho^s(n)}{n^{2-p/2} }
< \infty\qquad\mbox{where $s = 2(\alpha-2)/{\alpha}$}.
\]
Then the conclusion of Theorem \ref{Wasser} holds for any $1\le r \le2$.
\end{cor}

\begin{pf*}{Proof of Theorem \ref{Wasser}}
Notice first that (\ref{cond1Wasser}) implies in particular that
$\|\E_0 (S_n) \|_p = o(n^{2/p^2})$ (apply, e.g., item 2
of Lemma 37 in \cite{MP} to the sequence $(\|\E_0 (S_n) \|^2_p
)_{n \geq0}$). Now, since $p>2$, (\ref{cond1Wasser}) then entails that
(\ref{MWpalpha}) holds true. Therefore, by Theorem \ref{theoMW}, $D$
defined by (\ref{defdiffmart}) is in ${\mathbb L}^p$. In addition,
since $p>2$, (\ref{MWpalpha}) implies that $\sum_{n>0} n^{-3/2} \| \E_0
(S_n) \|_2 < \infty$ which is a sufficient condition for $n^{-1}\E
(S_n^2)$ to converge; see Theorem 1 in \cite{PU}.

Let now $M_n=\sum_{k=1}^{n}D \circ\theta^k$ and $R_n=S_n-M_n$.
According to the proof of Theorem 3.1 in \cite{DMRi} and to their
Remark 2.1, the theorem will follow if we can prove that
%
%
\begin{equation}
\label{neglrn} \|R_n\|_r=O\bigl(n^{(3-p)/2}\bigr)
\end{equation}
and also that
\[
\sum_{k \geq0} \frac{\|\E_{0}
(M_{2^k}^2 ) - \E
(M_{2^k}^2 )\|_{1 + \gamma}}{2^{k(2-p/2)}}< \infty\qquad\mbox{for a }
\gamma>0
\]
and
\[
\sum_{k \geq0}\frac{ \|\E_{0}
(M_{2^k}^2 ) - \E
(M_{2^k}^2 )\|_{p/2}}{2^{2k/p}}< \infty.
\]
Using (\ref{compdyanorm}) and the subadditivity of the sequence
$(\|\E_{-2n}
(M_{n}^2 ) - \E
(M_{n}^2 )\|_{q})_{n \geq1}$, for any $q \geq1$, we infer that
the latter conditions are equivalent to
%
%
\begin{eqnarray}
\label{condmartwasserdya} \sum_{n \geq1}
\frac{\|\E_{-2n}
(M_{n}^2 ) - \E
(M_{n}^2 )\|_{1 + \gamma}}{n^{3-p/2}}&<& \infty\qquad\mbox{for a
}\gamma>0
\quad\mbox{and}\nonumber\\[-8pt]\\[-8pt]
\sum
_{n \geq1}\frac{ \|\E_{-2n}
(M_{n}^2 ) - \E
(M_{n}^2 )\|_{p/2}}{n^{1+2/p}}&<& \infty.\nonumber
\end{eqnarray}
Using Proposition \ref{passMnSn} we infer that (\ref
{condmartwasserdya}) holds, provided that (\ref{condSnMDPnondya}) and
(\ref{Cond1carre}) do, and that
%
%
\begin{eqnarray}
\label{p1condmartwasser} \sum_{n \geq1}
\frac{\|\E_0(S_n) \|^2_{p}}{n^{1+4/p^2}}&<& \infty,\qquad \sum_{n \geq1}
\frac{\|\E_0(S_n) \|^2_{2(1
+ \gamma)}}{n^{1+(4-p)/(2+2\gamma)}}< \infty\quad\mbox{and}\nonumber\\[-8pt]\\[-8pt]
\sum_{n \geq1}
\frac{ \|\E_0(S_n) \|_{2}}{n^{(5-p)/2}}&<& \infty.\nonumber
\end{eqnarray}
Notice first that the third part of (\ref{p1condmartwasser}) holds,
provided that (\ref{cond2Wasser}) does [notice that the second part of
(\ref{cond2Wasser}), for $r>2$ implies the first part of (\ref
{cond2Wasser})], whereas the first part of (\ref{p1condmartwasser}) is
exactly condition (\ref{cond1Wasser}).
Notice now that for any $p \in\ ]2,3[$ and $\gamma$ small enough, $
(4-p)/(2 + 2 \gamma) \geq4/p^2$ and $p\geq2 + 2\gamma$. Therefore
the second part of (\ref{p1condmartwasser}) is implied by condition
(\ref{cond1Wasser}).

It remains to prove (\ref{neglrn}). By Lemma 2.7 of \cite{PU}, the
first part of (\ref{cond2Wasser}) implies that $\|\E_0(S_n)\|
_2=o(n^{(3-p)/2})$. Therefore by using Theorem \ref{theoMW}, we infer
that, since $p>2$, for any $r$ in $[1,2]$, $\|R_n \|_r \leq\|R_n\|
_2=o(n^{(3-p)/2})$ under the first part
of~(\ref{cond2Wasser}). Now, since $p>2$, for any $r$ in $]2,p]$, the
second part of (\ref{cond2Wasser})
implies that $\|R_n \|_r =O(n^{(3-p)/2})$ by Theorem \ref{theoMW}.
\end{pf*}

\section{Proof of the martingale approximation results} \label{sectproof}

In all the following lemmas, $p>1$ and $(X_n)_{n\in\mathbb{Z}}$ is an
adapted stationary sequence in $\LL^p({\mathcal H})$ in the sense of
Notation \ref{notation1}.

%
%
\begin{lem}\label{lem}
We have $\sum_{k\ge0} (k+1)^{-1}\|{\mathcal P}_0(X_k)
\|_{p,{\mathcal H}}<\infty$.
\end{lem}
\begin{pf}
We first prove the case $p\ge2$. By H\"older's
inequality, we have
\begin{eqnarray*}
\biggl(\sum_{k\ge0}\frac{\|{\mathcal P}_0(X_k)\|
_{p,{\mathcal H}}}{k+1}
\biggr)^p &\ll&\sum_{k\ge0}\bigl\|{\mathcal
P}_{-k}(X_0) \bigr\|_{p,{\mathcal H}}^p \ll\biggl\|
\biggl(\sum_{k\ge
0}\bigl|{\mathcal P}_{-k}(X_0)\bigr|_{{\mathcal H}}^p
\biggr)^{1/p} \biggr\|_{p}^p
\\
&\ll& \biggl\| \biggl(\sum_{k\ge0}\bigl|{\mathcal
P}_{-k}(X_0)\bigr|_{{\mathcal
H}}^2
\biggr)^{1/2} \biggr\|_{p}^p \ll\|X_0
\|_{p,{\mathcal H}}^p,
\end{eqnarray*}
where we used $\|\cdot\|_{\ell^p}\le\|\cdot\|_{\ell^2} $
and Burkholder's inequality for ${\mathcal H}$-valued martingales; see
\cite{BU}.

Let prove the case $1<p<2$.
By H\"older's inequality,
\begin{eqnarray*}
\biggl(\sum_{k\ge0}\frac{\|{\mathcal P}_0(X_k)\|
_{p,{\mathcal H}}}{k+1}
\biggr)^p &\ll&\sum_{k\ge0}
\frac{\|{\mathcal
P}_{-k}(X_0)
\|_{p,{\mathcal H}}^p}{(k+1)^{p/2}}= \E\biggl(\sum_{k\ge0}
\frac
{|{\mathcal P}_{-k}(X_0)|_{{\mathcal H}}^p} {
(k+1)^{p/2}} \biggr)
\\
&\ll& \biggl\| \biggl(\sum_{k\ge0}\bigl|{\mathcal
P}_{-k}(X_0)\bigr|_{{\mathcal
H}}^2
\biggr)^{1/2} \biggr\|_p^p \ll\|X_0
\|_{p, {\mathcal H}}^p,
\end{eqnarray*}
where we used again H\"older's inequality and Burkholder's inequality for
${\mathcal H}$-valued martingales.
\end{pf}

%
%
\begin{lem}\label{lemme0}
Assume that
%
%
\begin{equation}
\label{cond} \sum_{n\ge1} \sum
_{k\ge0} \frac{\|{\mathcal P}_{0}(S_n\circ\theta^{k-1})
\|_{p,{\mathcal H}}}{(n+k)^2}<\infty.
\end{equation}
Then $\sum_{n\ge0}|\sum_{k\ge n} \frac{{\mathcal
P}_0(X_k)}{k+1}|_{{\mathcal H}}$ converges in
$\LL^p$ and a.s. Moreover for any integer $m \geq0$,
%
%
\begin{equation}
\label{est1} \biggl\|\sum_{n\ge m}\sum
_{k\ge n} \frac{{\mathcal P}_0(X_k)}{k+1} \biggr\|_{p,{\mathcal H}} \le\sum
_{k\ge m}\sum_{n\ge1}
\frac{\|{\mathcal P}_{0}(S_n\circ
\theta^{k-1})
\|_{p,{\mathcal H}}}{(n+k)^2}.
\end{equation}
\end{lem}
\begin{pf}
By assumption, the series
\[
\sum_{k\ge0} \biggl|\sum_{n\ge1}
\frac{\sum_{l=0}^{n-1}{\mathcal
P}_0(X_{l+k})
}{(n+k)(n+k+1)} \biggr|_{{\mathcal H}}
\]
converges a.s. and in $\LL^p$. On the other hand, using Lemma \ref
{lem} to invert the order of summation, we have
\begin{eqnarray*}
\sum_{l\ge k} \frac{{\mathcal P}_0(X_{l})}{l+1} &=&\sum
_{l\ge0}\frac{{\mathcal
P}_0(X_{k+l})}{k+l+1} = \sum_{l\ge0}
\sum_{n\ge l+1} \frac{{\mathcal P}_0(X_{k+l})}{(k+n)(k+n+1)}\\
&=& \sum
_{n\ge1}\frac
{\sum_{l=0}^{n-1}{\mathcal P}_0(X_{l+k})
}{(n+k)(n+k+1)},
\end{eqnarray*}
which gives the desired convergence.
\end{pf}

%
%
\begin{lem}\label{lemmedalibor}
For every integer $r \geq0$,
%
%
\begin{equation}
\label{est4} \sum_{k\ge r}\sum
_{m\ge1} \frac{\|{\mathcal P}_{-k}(S_m)
\|_{p, {\mathcal H}}}{(m+k)^2} \ll\sum
_{k \ge r+1}\frac{\|\E_{-r}(S_k)\|_{p, {\mathcal
H}}}{k^{1+1/p''}}.
\end{equation}
\end{lem}
\begin{pf}
Let $m$ be a positive integer. Assume first that $p\ge2$.
By H\"older's inequality
and using that $\|\cdot\|_{\ell^p}\le\|\cdot\|_{\ell^2}$, we have
\begin{eqnarray*}
\sum_{k\ge r} \frac{\|{\mathcal P}_{-k}(S_m)
\|_{p, {\mathcal H}}}{(m+k)^2} &\ll&(m+r)^{-(1+1/p)}
\biggl(\sum_{k\ge r} \bigl\|{\mathcal P}_{-k}(S_{m})
\bigr\|_{p, {\mathcal H}}^p \biggr)^{1/p}
\\
&\ll&(m+r)^{-(1+1/p)} \biggl( \E\biggl(\sum_{k\ge r}
\bigl|{\mathcal P}_{-k}(S_{m})\bigr|_{{\mathcal H}}^2
\biggr)^{p/2} \biggr)^{1/p} \\
&\ll&\frac{\|
\E_{-r}(S_m)\|_{p, {\mathcal H}}}{(m+r)^{1+1/p}},
\end{eqnarray*}
where we used Burkholder's inequality for ${\mathcal H}$-valued
martingales (see \cite{BU}), in the last step.

Assume now that $1<p< 2$. We use H\"older's inequality twice and once
again Burkholder's inequality for ${\mathcal H}$-valued martingales in
the last step, to obtain
\begin{eqnarray*}
\sum_{k\ge r} \frac{\|{\mathcal P}_{-k}(S_m)
\|_{p, {\mathcal H}}}{(m+k)^2} &\ll&\frac{1}{(m+r)^{1/p}}
\biggl( \sum_{k\ge r} \frac{\|{\mathcal
P}_{-k}(S_{m})\|_{p, {\mathcal H}}^p}{(m+k)^p}
\biggr)^{1/p}
\\
&\ll&\frac{1}{(m+r)^{1/p}} \biggl( \frac{1}{(m+r)^{3p/2-1}}\E\biggl(
\sum
_{k\ge r} \bigl|{\mathcal P}_{-k}(S_{m})\bigr|_{{\mathcal H}}^2
\biggr)^{p/2} \biggr)^{1/p}\\
&\ll&\frac{\|
\E_{-r}(S_m)\|_{p, {\mathcal H}}}{(m+r)^{3/2}}.
\end{eqnarray*}
From the above computations, we then derive that
\begin{eqnarray*}
\sum_{k\ge r} \sum_{m\ge1}
\frac{\|{\mathcal P}_{-k}(S_m)
\|_{p, {\mathcal H}}}{(m+k)^2} & \ll& \sum_{m\ge1} \frac{\|\E
_{-r}(S_m)\|_{p, {\mathcal H}}}{(m+r)^{1+1/p''}}
\\
& \ll&\frac{1}{(r+1)^{1/p''}}\max_{1 \leq m \leq r} \bigl\|\E_{-r}(S_m)
\bigr\| _{p, {\mathcal H}}\\
&&{}+ \sum_{m\ge r+1} \frac{\|\E_{-r}(S_m)\|_{p,
{\mathcal H}}}{m^{1+1/p''}}.
\end{eqnarray*}
The lemma then follows by using Lemma \ref{SUB} with $\gamma=1/p''$
and $\ell=r$.
\end{pf}

%
%
\begin{lem}\label{lemme}
For every $r\ge0$,
%
%
\begin{eqnarray}
\label{decomp} X_0&=& \sum_{k=0}^r
\sum_{l\ge k}\frac{{\mathcal P}_0(X_l)}{l+1} -\sum
_{k=0}^r\sum_{l\ge k}
\frac{\E_0(X_{l+1})-\E_{-1}(X_l)}{l+1}\nonumber\\[-8pt]\\[-8pt]
&&{}+(r+1) \sum_{l\ge r}\frac{\E
_0(X_{l+1})}{(l+1)(l+2)}.\nonumber
\end{eqnarray}
In particular, if we assume (\ref{MWpalpha}), letting $r\to\infty$,
we have
\[
X_0= \sum_{k\ge0}\sum
_{l\ge k}\frac{{\mathcal P}_0(X_l)}{l+1} -\sum_{k\ge0}
\sum_{l\ge k}\frac{\E_0(X_{l+1})-\E_{-1}(X_l)}{l+1}.
\]
\end{lem}
\begin{pf}
Let $m\ge k\ge0$. We have
\[
\sum_{l=k}^m\frac{{\mathcal P}_0(X_{l})}{l+1}=
\frac{\E_0(X_k)}{k+1} -\frac{\E_0(X_{m+1})}{m+2}+\sum_{l=k}^m
\frac{\E
_0(X_{l+1})}{l+2}-\sum_{l=k}^m
\frac{\E_{-1}(X_{l})}{l+1}.
\]
Hence
\begin{eqnarray*}
\sum_{l=k}^m\frac{{\mathcal P}_0 (X_{l})}{l+1}&=&
\frac{\E_0(X_k)}{k+1} -\frac{\E_0(X_{m+1})}{m+2}+\sum_{l=k}^m
\frac{\E_0(X_{l+1})-\E_{-1}(X_{l})} {
l+1}\\
&&{}-\sum_{l=k}^m
\frac{\E_{0}(X_{l+1})}{(l+1)(l+2)}.
\end{eqnarray*}
Notice that $m^{-1}\|\E_0(X_m)\|_{p, {\mathcal H}}
\to0$ and that $\sum_{l\ge
0}\frac{\|\E_{0}(X_{l+1})\|_{p,{\mathcal H}}}{(l+1)(l+2)}<\infty$.
Hence, using Lemma \ref{lem}, we may and do let
$m\to\infty$, to obtain
\[
\sum_{l\ge k}\frac{{\mathcal P}_0(X_{l})}{l+1}= \frac{\E_0(X_k)}{k+1}+
\sum_{l\ge k}\frac{\E_0(X_{l+1})-\E_{-1}(X_{l})} {
l+1}-\sum
_{l\ge k}\frac{\E_0(X_{l+1})}{(l+1)(l+2)}.
\]
Let $r\ge0$. We then deduce that
\begin{eqnarray*}
\sum_{k=0}^r\sum
_{l\ge k}\frac{{\mathcal P}_0(X_{l})}{l+1}&=&\sum_{k=0}^r
\frac{\E_0(X_k)}{k+1}+\sum_{k=0}^r\sum
_{l\ge k}\frac{\E_0(X_{l+1})-
\E_{-1}(X_{l})} {
l+1}\\
&&{}-\sum
_{k=0}^r\sum_{l\ge k}
\frac{\E_0(X_{l+1})}{(l+1)(l+2)}.
\end{eqnarray*}
Hence, interverting the order of summation in the last term,
\begin{eqnarray*}
\sum_{k=0}^r\sum
_{l\ge k}\frac{{\mathcal P}_0(X_{l})}{l+1}&=& X_0+\sum
_{k=0}^r\sum_{l\ge k}
\frac{\E_0(X_{l+1})-\E_{-1}(X_{l})} {
l+1}\\
&&{}- (r+1)\sum_{l\ge r}
\frac{\E_0(X_{l+1})}{(l+1)(l+2)}.
\end{eqnarray*}
Assume (\ref{MWpalpha}). In view of Lemmas \ref{lemme0} and \ref
{lemmedalibor},
we see that the series on the left converges in $\LL^p(\H)$. On an
other hand, Lemma \ref{SUB} (with $\gamma=1$) implies that $n^{-1} \|
\E_0(S_n)\|_{p,\H}\to0$. Therefore by Abel summation,
\[
\biggl\|(r+1)\sum_{l\ge r}\frac{\E_0(X_{l+1})} {
(l+1)(l+2)}\biggr\|_{p,\H}\to0,
\]
when $r\to\infty$.
\end{pf}

\subsection{\texorpdfstring{Proof of Theorem \protect\ref{theoMW}}{Proof of Theorem 2.3}}
The first assertion comes from Lemma \ref{lemme0} combined with Lemma
\ref{lemmedalibor}. Now, by Lemma \ref{lemme}, we have
\[
X_1= D\circ\theta-\sum_{k\ge0}\sum
_{l\ge k+1}\frac{\E_1(X_{l+1})-\E
_{0}(X_l)}{l}.
\]
Hence, using that $\E_1(X_{l+1})=\E_{0}(X_l)\circ\theta$, we obtain
that for any positive integer~$n$,
\begin{eqnarray*}
S_n-M_n&=&-\sum_{k\ge0}\sum
_{l\ge k+1}\frac{
\E_0(X_{l})\circ\theta^{n}-\E_{0}(X_l)}{l} \\
&=&-\sum
_{k\ge0}\sum_{l\ge k+1}
\frac{\E_{n}(X_{l+n})-\E
_{0}(X_l)}{l}.
\end{eqnarray*}
Let $N$ be a positive integer, fixed for the moment. Then writing
%
%
\begin{equation}
\label{Vn2} V_{n,N}=\sum_{k=0}^{N-1}
\sum_{l\ge k+1}\frac{\E_{n}(X_{l+n})-\E
_{0}(X_{l+n})}{l}
\end{equation}
and
%
%
\begin{equation}
\label{Wn2} W_{n,N} =\sum_{k\ge N}
\sum_{l\ge k+1}\frac{\E_{n}(X_{l+n})-\E
_{0}(X_{l+n})}{l},
\end{equation}
we obtain
%
%
\begin{eqnarray}
\label{dec} S_n-M_n-\E_{0}(S_n)&=&-
\sum_{k\ge0}\sum_{l\ge k+1}
\frac{\E_{n}
(X_{l+n})-\E_{0}(X_{l+n})}{l}\nonumber\\[-8pt]\\[-8pt]
&=&-(V_{n,N}+W_{n,N}).\nonumber
\end{eqnarray}
%

We first deal with $V_{n,N}$. We have
%
%
\begin{eqnarray}\label{Vn}\quad
V_{n,N}&=&\sum_{l=1}^{N}
\bigl( \E_{n}(X_{l+n})-\E_{0}(X_{l+n})
\bigr) + N\sum_{l\ge N+1} \frac{\E_{n}(X_{l+n})-\E_{0}(X_{l+n})}{l}
\nonumber\\[-8pt]\\[-8pt]
&=&\E_{0}(S_{N})\circ\theta^n-
\E_{0}\bigl(S_{N}\circ\theta^n\bigr) +N\sum
_{l\ge N+1} \frac{\E_{n}(X_{l}\circ\theta^n)-\E_{0}
(X_{l}\circ\theta^n)}{l}.\nonumber
\end{eqnarray}
Let $j\in\{0,n\}$. By (\ref{MWpalpha}) and Lemma \ref{SUB} with
$\gamma=1$,
%
%
\begin{equation}
\label{NN} \frac{\|\E_0(S_{N})\|_{p,{\mathcal H}}}{N}\ll\sum_{l\ge N}
\frac{\|\E_{0}
(S_{l})\|_{p,{\mathcal H}}}{l^2}=o(1).
\end{equation}
Using Abel summation
we have, for every $s\ge N+1$,
\begin{eqnarray*}
\sum_{l= N+1}^s
\frac{\E_{j}(X_{l}\circ\theta^n)}{l}&=& \sum_{l= N+1}^s
\frac{\E_{j}(S_{l}\circ\theta^n-S_{l-1}
\circ\theta^n)}{l}
\\
&=&-\frac{\E_j(S_{N}\circ\theta^n)}{N+1}+\frac{\E_j(S_{s}\circ
\theta^n)}{s+1}+\sum_{l= N+1}^s
\frac{\E_{j}
(S_{l}\circ\theta^n)}{l(l+1)}.
\end{eqnarray*}
Letting $s\to\infty$, it follows from (\ref{NN}) that
%
%
\begin{equation}
\label{calcul} \sum_{l\ge N+1} \frac{\E_{j}(X_{l}\circ\theta^n)}{l} = -
\frac{\E_j(S_{N}\circ\theta^n)}{N+1}+\sum_{l\ge N+1}\frac{\E_{j}
(S_{l}\circ\theta^n)}{l(l+1)}.
\end{equation}
Hence, starting from (\ref{Vn}) and considering (\ref{calcul}) and
(\ref{NN}), we derive that
%
%
\begin{eqnarray}
\label{boundVn}\quad
\| V_{n,N}\|_{p,{\mathcal H}} &\le& 2
\frac{\|\E_0(S_{N})\|
_{p,{\mathcal H}}}{N}+ N \sum_{l\ge N+1}\frac{\|\E_{n}
(S_{l}\circ\theta^n)-\E_{0}
(S_{l}\circ\theta^n)\|_{p,{\mathcal H}}}{l(l+1)}
\nonumber\\[-8pt]\\[-8pt]
&\ll& N \sum_{l\ge N}\frac{\|\E_{0}
(S_{l})\|_{p,{\mathcal H}}}{l^2}.\nonumber
\end{eqnarray}
It remains to deal with $W_{n,N}$. Since $\E_{0}(W_{n,N})=0$, we have
\[
W_{n,N}=\sum_{r=1}^{n} {\mathcal P}_r(W_{n,N}).
\]
Using that ${\mathcal P}_r$ defines a continuous operator on
$\LL^p(\H)$ and that the series in (\ref{Wn2}) converges in $\LL
^p(\H)$,
we infer that
%
%
\begin{equation}
\label{Wnbis} W_{n,N} =\sum_{r=1}^{n}
\sum_{k\ge N}\sum_{l\ge k+1}
\frac{\E_{r}(X_{l+n})-\E
_{r-1}(X_{l+n})}{l}.
\end{equation}
But, by Burkholder's inequality for ${\mathcal H}$-valued martingales
(see \cite{BU}),
%
%
\begin{equation}
\label{Wn} \|W_{n,N}\|_{p,{\mathcal H}}^{p'}\ll\sum
_{r=1}^{n} \bigl\|{\mathcal P}_r(W_{n,N})
\bigr\|_{p,{\mathcal H}}^{p'}.
\end{equation}
Notice that for any $r \in\{1,\ldots, n \}$,
\[
{\mathcal P}_r(W_{n,N}) = \biggl(\sum
_{k\ge N}\sum_{l\ge1}
\frac
{{\mathcal
P}_0(X_{l+k+n-r})}{l+k} \biggr)\circ\theta^{r}.
\]
Now, using Lemma \ref{lem},
\begin{eqnarray*}
\sum_{l\ge1}\frac{{\mathcal
P}_0(X_{l+k+n-r})}{l+k} &=& \sum
_{l\ge1}{\mathcal P}_0(X_{l+k+n-r})\sum
_{m\ge l}\frac{1}{(m+k)(m+k+1)} \\
&=& \sum
_{m\ge
1}\frac{{\mathcal P}_0(S_m\circ\theta^{k+n-r})} {
(m+k)(m+k+1)}.
\end{eqnarray*}
Therefore,
%
%
\begin{equation}
\label{difmart} \biggl|\sum_{k\ge N}\sum
_{l\ge1}\frac{{\mathcal
P}_0(X_{l+k+n-r})}{l+k} \biggr|\le\sum
_{m\ge1} \sum_{k\ge N}
\frac{|{\mathcal P}_0(S_m\circ\theta^{k+n-r})|} {
(m+k)^2}.
\end{equation}
Hence, with $s=n-r$,
\[
\|W_{n,N}\|_{p,{\mathcal H}}\ll n^{1/p'} \max
_{0\le s \le n-1} \sum_{k\ge N+s}\sum
_{m\ge1} \frac{\|{\mathcal P}_{-k}(S_m)\|_{p,{\mathcal H}}}{(m+k-s)^2}.
\]
Now we take $N=u_n \ge n$. We then infer that
%
%
\begin{equation}
\label{boundWn} \|W_{n,u_n}\|_{p,{\mathcal H}}\ll n^{1/p'} \sum
_{k\ge u_n}\sum_{m\ge1}
\frac{\|{\mathcal P}_{-k}(S_m)\|_{p,{\mathcal H}}}{(m+k)^2}.
\end{equation}
Hence using (\ref{dec}), (\ref{boundVn}) with $N=u_n$ and (\ref
{boundWn}), we get that
%
%
\begin{eqnarray}
\label{1decth26}
\|S_n-M_n\|_{p, {\mathcal H}} &\ll&\bigl\|
\E_0(S_n)\bigr\|_{p, {\mathcal H}} + u_n\sum
_{m\ge
u_n}\frac{\|\E_0(S_m)\|_{p, {\mathcal H}}}{m^2}\nonumber\\[-8pt]\\[-8pt]
&&{} +n^{1/p'}\sum
_{k\ge u_n}\sum_{m\ge1}
\frac{\|{\mathcal P}_{-k}(S_m)\|_{p, {\mathcal H}}}{(m+k)^2}.\nonumber
\end{eqnarray}

Next using Lemma \ref{SUB} with $\gamma=1$, we derive that
%
%
\begin{equation}
\label{applisub} \bigl\|\E_0(S_n)\bigr\|_{p, {\mathcal H}} \leq\max
_{1 \leq k \leq u_n }\bigl\|\E_0(S_k)\bigr\|_{p, {\mathcal H}}
\ll u_n\sum_{m\ge
u_n}\frac{\|\E_0(S_m)\|_{p, {\mathcal H}}}{m^2}.
\end{equation}
Starting from (\ref{1decth26}) with $u_n = [n^q]$ and taking into
account (\ref{applisub}) and Lemma \ref{lemmedalibor}, Theorem \ref
{theoMW} follows.

\subsection{\texorpdfstring{Proof of Theorem \protect\ref{quenched}}{Proof of Theorem 2.7}}

Part of the proof relies on a new ergodic theorem with rate. Hence we
first recall some facts from ergodic
theory and state our ergodic theorem, while we give its proof in
Appendix \ref{sectergodic}.

Let $T$ be a Dunford--Schwartz operator on $\Omega$;
that is, $T$ is a contraction of $\LL^1$ and $\LL^\infty$. Let ${\mathbf
T}$ be the linear modulus of $T$; see, for example, Theorem 1.1,
Chapter 4 of \cite{Krengel}. Recall
that ${\mathbf T}$ is a positive Dunford--Schwartz operator such that
$|Tf|\le{\mathbf T}|f|$ for every $f\in\LL^1$ and
$|Tf|^p\le{\mathbf T}(|f|^p)$ for every $f\in\LL^p$.

We will make use, for $p\ge1$, of the weak $\LL^p$-spaces
\[
\LL^{p,w}:= \Bigl\{f\in\LL^0\dvtx \sup_{\lambda>0}
\lambda^p\P\bigl\{ |f| \ge\lambda\bigr\}<\infty\Bigr\},
\]
where $\LL^0$ is the space of all ${\mathcal A}-{\mathcal B} ( \mathbb
R)$ measurable functions.

Recall that when $p>1$, there exists a norm $\|\cdot\|_{p,w}$ on $\LL^{p,w}$
that\break  makes $\LL^{p,w}$ a Banach space and which is equivalent to the
``pseudo''-norm
$(\sup_{\lambda>0} \lambda^p\P\{ |f| \ge\lambda\})^{1/p}$.

We define, for every $l\ge0$, a maximal operator as follows. For any
nonnegative function $h\in\LL^1$, let
\[
\M_l(h)=\sup_{n\ge1} \frac{h+{\mathbf T}^{2^l}h +\cdots+({\mathbf
T}^{2^l})^{n-1}h}{n}.
\]
By the Dunford--Schwartz (or Hopf) ergodic theorem (see, e.g., Krengel
\cite{Krengel}, Lemma 6.1, page 51, and Corollary 3.8, page 131),
\[
\sup_{\lambda>0} \lambda\P\bigl\{\M_l(h)\ge\lambda
\bigr\}\le\|h\|_1.
\]
In particular, for every $p> 1$, there exists $C_p>0$ such that, for
every $f\in\LL^p$,
%
%
\begin{equation}
\label{inemax} \bigl\|\bigl(\M_l\bigl(|f|^p\bigr)
\bigr)^{1/p}\bigr\|_{p,w} \le C_p \|f\|_p.
\end{equation}

Let $\B$ be a Banach space with norm \mbox{$|\cdot|_\B$}. For every $p\ge
1$, we denote by
$\LL^p(\B)$ the Bochner space
$\{f\dvtx\Omega\to\B, |f|_\B\in\LL^p\}$. When $T$ is induced by a
measurable transformation $\theta$ preserving $\P$, $\M_l(|f|_\B)$
is well defined
for every $f\in\LL^1(\B)$. We prove the following, where
$U_n(f)=f+\cdots+ T^{n-1}f$.

%
%
\begin{prop}\label{prop}
Let $T$ be a Dunford--Schwartz operator on $(\Omega,{\mathcal A},\P)$
and $f\in\LL^1$.
We have
\[
\max_{1\le n \le2^r} \bigl|U_n(f)\bigr|\le2^{r/p} \sum
_{k=0}^r \frac{
[\M_k(|U_{2^k}(f)|^p) ]^{1/p}}{2^{k/p}}.
\]
When $T$ is induced by a measure preserving transformation $\theta$,
and $\B$ is a Banach space, the result holds also for $f\in\LL^1(\B)$,
replacing $|\cdot|$ with $|\cdot|_\B$.
\end{prop}
\begin{pf}
The proof follows from the following lemma,
using that
$U_{2^km}(f)-U_{2^k(m-1)}(f)=T^{2^k(m-1)}f+\cdots+
T^{2^km-1}f=(T^{2^k})^{(m-1)} U_{2^k}(f)$.
\end{pf}

%
%
\begin{lem}\label{Wu}
Let $(a_n)$ be a sequence in a Banach space $\B$ with norm $|\cdot
|_\B$. Write $s_n=a_1+ \cdots+ a_n$ and $s_0=0$. Let $p\ge1$. For
every $r\ge0$, we have
%
%
\begin{equation}
\max_{1\le n \le2^r} |s_n|_\B\le\sum
_{k=0}^r \Biggl(\sum_{m=1}^{2^{r-k}}
|s_{2^km}-s_{2^{k}(m-1)}|_\B^p
\Biggr)^{1/p}.
\end{equation}
\end{lem}
\begin{pf}
We make the proof by induction on $r\ge0$. The result
is obvious for $r=0$. Let $1 \le n \le2^r $. We have
$|s_{2n-1}|_\B\le|s_{2n-2}|_\B+|a_{2n-1}|_\B$. Hence,
writing $\widetilde a_n =a_{2n-1}+a_{2n}$ and $\widetilde s_n =
\sum_{k=1}^n \widetilde a_k=s_{2n}$, we get that
\[
\max_{1\le l\le2^{r+1}}|s_l|_\B\le\max
_{1\le n \le2^r}|\widetilde s_n|_\B+ \Biggl(\sum
_{n=1}^{2^r}|a_{2n-1}|_\B^p
\Biggr)^{1/p},
\]
and the result follows.
\end{pf}

%
%
\begin{theo}\label{theoergo}
Let $T$ be a Dunford--Schwartz operator on $(\Omega,{\mathcal A},\P
)$. Let
\mbox{$f\in\LL^p$}, $p>1$. Let $\psi$ be a positive nondecreasing\vadjust{\goodbreak}
function, such that there exists $C>1$
such that $\psi(2x )\le C\psi(x)$, for every $x\ge1$. Assume that
%
%
\begin{equation}
\label{MW1} \sum_n \frac{\|f+\cdots+ T^{n-1}f\|_p} {\psi(n) n^{1+1/p}}<
\infty.
\end{equation}
Then $\sup_{n\ge1}\frac{|f+\cdots+ T^{n-1}f|}{\psi(n)n^{1/p}}\in
\LL^{p,w}$ and
$ \frac{|f+\cdots+T^{n-1}f|}{\psi(n) n^{1/p} }\to0$
$\P$-a.s.\vspace*{1pt}

If $T$ is induced by a measure-preserving transformation, and $(\B,
|\cdot|_\B)$ is a Banach space, the result holds
with $|\cdot|_\B$ instead of $|\cdot|$ for every $f\in\LL^p(\B)$
such that $\sum_n \frac{\| |f+\cdots+ T^{n-1}f|_\B\|_{p}} {
\psi(n) n^{1+1/p}}< \infty$.
\end{theo}
%
%

We turn now to the proof of Theorem \ref{quenched}. It will follow
from the next two propositions. Notice that the second one is a version
of Corollary
22 of Merlev\`ede and Peligrad \cite{MP} under $\E_0$.

%
%
\begin{prop}\label{propquenched1}
Assume (\ref{MW}). Then
$\E_0[(S_n-M_n-\E_0(S_n))^2] =o(n)$ $\P$-a.s. and $\E
_0(S_n)=o(\sqrt n)$ $\P$-a.s. In
particular,
\[
\E_0\bigl[(S_n-M_n)^2\bigr]
=o(n)\qquad \mbox{$\P$-a.s.}
\]
\end{prop}

%
%
\begin{prop}\label{propquenched2}
Assume (\ref{MW}) and that $\E_0(S_n^2)=o(n)$ $\P$-a.s. Then
%
%
\begin{equation}
\E_0\Bigl(\max_{1\le k\le n} S_k^2
\Bigr)=o(n)\qquad \mbox{$\P$-a.s.}
\end{equation}
\end{prop}

Before proving the above propositions, we indicate how they lead to
Theorem~\ref{quenched}. Using Proposition \ref{propquenched1}, we
apply Proposition \ref{propquenched2} with $S_n-M_n$ in place of
$S_n$. This proves (\ref{resquen}). Now the convergence (\ref{FCLT})
follows from (\ref{resquen}) together with the quenched weak
invariance principle for martingales; see, for instance, Derriennic and
Lin \cite{DL} for the ergodic case. To be more precise, if we define
$D_k=D\circ\theta^k$ and $\widetilde W_n$ by $\widetilde
W_n(t)=n^{-1/2}(M_{[nt]}+(nt-[nt])D_{[nt]+1})$, then (\ref{FCLT})
holds with $\widetilde W_n$ in place of $ W_n$, and $ \eta= \E( D^2
|{\mathcal I})$.
To end the proof, we first notice that by Theorem 1 of Peligrad and
Utev \cite{PU}, $\E( D^2 |{\mathcal I}) = \lim_{n \rightarrow\infty}
n^{-1} \E( S_n^2|{\mathcal I})$ in $\LL^1$. It remains to prove that
$\E( D^2 |{\mathcal I}) = \lim_{n \rightarrow\infty} n^{-1}
\E_0 ( S_n^2)$ in $\LL^1$. But, by (\ref{MW}) and (\ref{rate}), $\|
S_n^2-M_n^2\|_1= o(n)$. Hence it suffices to prove that $\E( D^2
|{\mathcal I}) = \lim_{n \rightarrow\infty} n^{-1}
\E_0 ( M_n^2)$ in $\LL^1$.

With this aim, we will make use of the operator $Q$ defined
by
\[
QZ=\E_0(Z\circ\theta) \qquad\forall Z\in\mathbb{L}^1.
\]
The operator $Q$ is Markovian and hence is a Dunford--Schwartz operator.
Notice that $Q^nZ=\E_0(Z\circ\theta^n)$. Moreover, by Lemma 7.1 in
\cite{DMP}, if $Z$ is additionnally assumed to be in ${\mathcal
F}_{\infty}$,
%
%
\begin{equation}
\label{ergthmQ} \mbox{$\bigl(QZ+\cdots+Q^nZ\bigr)/n$ converges $
\P$-a.s. and in $\LL^1$ to $\E(Z|{\mathcal I})$}.
\end{equation}

To conclude we take $Z=D^2$ and we notice that, by orthogonality,
$\E_0(M_n^2)=Q(D^2)+ \cdots+ Q^{n}(D^2)$.

It remains to prove Propositions \ref{propquenched1} and \ref{propquenched2}.

\begin{pf*}{Proof of Proposition \ref{propquenched1}}
The fact that $\E_0(S_n)=o(\sqrt n)$ $\P$-a.s. under (\ref{MW}) comes
directly from an application of Theorem \ref{theoergo} with $T=Q$. We
prove now that under (\ref{MW}), the following convergence holds: $\E
_0[(S_n-M_n-\E_0(S_n))^2] =o(n)$ $\P$-a.s.

Let $N$ be a positive integer fixed for the moment. By
(\ref{dec}), we have
%
%
\begin{equation}
\label{eqt} S_n-M_n-\E_{0}(S_n)=-(V_{n,N}+W_{n,N}),
\end{equation}
where $V_{n,N}$ and $W_{n,N}$ are given, respectively, by (\ref{Vn2})
and (\ref{Wn2}).

Let $\varphi_N:=\E_0(S_N)$ and $\psi_N=\sum_{l\ge N+1} \frac
{\varphi_l}{l(l+1)}$, where $\psi_N$ is well defined in $\mathbb
{L}^2$, by (\ref{MW}).

Then, by (\ref{Vn}) and (\ref{calcul}),
\[
|V_{n,N}|\ll\bigl|\varphi_N\circ\theta^n \bigr|
+\bigl|Q^n\varphi_N\bigr|+\bigl|\psi_N\circ
\theta^n \bigr|+\bigl|Q^n\psi_N\bigr|.
\]
Hence, by using (\ref{ergthmQ}),
\[
\E_0\bigl(V_{n,N}^2\bigr)\ll Q^n
\bigl(\varphi_N^2\bigr) + Q^n\bigl(
\psi_N^2\bigr) = o(n) \qquad\mbox{$\P$-a.s.}
\]
Then, using that $\E_0(S_n)= o(\sqrt n)$ $\P$-a.s. and (\ref{eqt}),
we obtain
\[
\limsup_n \frac{\E_0((S_n-M_n)^2)}n\le\limsup
_n \frac{\E
_0(W_{n,N}^2)}n.
\]

It remains to deal with $W_{n,N}$. Recall that by (\ref{Wnbis}),
\[
W_{n,N}=\sum_{r=1}^{n} {
\mathcal P}_r(W_{n,N})=\sum_{r=1}^{n}
\biggl(\sum_{k\ge N}\sum_{l\ge1}
\frac{{\mathcal
P}_0(X_{l+k+n-r})}{l+k} \biggr)\circ\theta^{r}.
\]
Hence, by orthogonality,
\[
\E_0\bigl(W_{n,N}^2\bigr)= \sum
_{r=1}^{n} \E_0\bigl({\mathcal
P}_r(W_{n,N})^2\bigr) =\sum
_{r=1}^{n}Q^r \biggl(\sum
_{k\ge N}\sum_{l\ge1}
\frac{{\mathcal
P}_0(X_{l+k+n-r})}{l+k} \biggr)^2.
\]
But, using (\ref{difmart}) and Cauchy--Schwarz's inequality, we have
\[
\Biggl|\sum_{k\ge N}\sum_{l\ge1}
\frac{{\mathcal
P}_0(X_{l+k+n-r})}{l+k} \Biggr| \ll\sum_{m\ge1}
\frac
{1}{(m+N)^{3/2}} \biggl(\sum_{k\ge0} \bigl|{\mathcal
P}_{-k}(S_m )\bigr|^2\circ\theta^k
\biggr)^{1/2}.
\]
Let now $g_N:=\sum_{m\ge1} \frac{1}{(m+N)^{3/2}}(\sum_{k\ge0}
|{\mathcal P}_{-k}(S_m )|^2\circ\theta^k)^{1/2}$.
Then $g_N$ is in $\mathbb{L}^2$ and
\[
\|g_N\|_2 \le\sum_{m\ge1}
\frac{\|\E_0(S_m)\|
_2}{(m+N)^{3/2}}<\infty.
\]
In particular, $\|g_N\|_2\to0$, as $N\to\infty$. So, finally, by
using (\ref{ergthmQ}), we get that
\[
\frac{\E_0(W_{n,N}^2)}{n}\ll\frac{\sum_{r=1}^n Q^r(g_N^2)}{n} \underset
{n\rightarrow+\infty} {
\longrightarrow} \E\bigl(g_N^2|{\mathcal I}\bigr)\qquad
\mbox{$\P$-a.s.}
\]
Since $\|\E(g_N^2|{\mathcal I})\|_1\le\|g_N^2\|_1\to0$, there exists
a sub-sequence $(N_j)$ such that
$\E(g_{N_j}^2|{\mathcal I})\rightarrow0$ $\P$-a.s. as $j\to\infty
$, and the result follows.
\end{pf*}

To prove Proposition \ref{propquenched2}, we will make use of the
following maximal inequality from
Merlev\`ede and Peligrad \cite{MP}. They did not state the result exactly
in that context,
but it may be proved exactly the same way, applying Doob's maximal inequality
conditionally, so the proof is omitted.

%
%
\begin{prop}\label{propMP}
Let $(X_n)_{n\in\mathbb{Z}}$ be a stationary sequence in $\LL^2$ in
the sense of
Notation \ref{notation1} and adapted to the filtration
$(\F_n)$. We have
%
%
\begin{eqnarray}
\label{MP}
&&
\Bigl(\E_0 \Bigl(\max_{1\le i\le2^r }
|S_i|^2 \Bigr) \Bigr)^{1/2}\nonumber \\
&&\qquad \le 2 \bigl(
\E_0\bigl(S_{2^r}^2\bigr) \bigr)^{1/2}+
2 \sum_{l=0}^{r-1} \Biggl( \sum
_{k=1}^{2^{r-l}-1} \E_0 \bigl( \bigl(\E
_{k2^l}(S_{(k+1)2^l}) -S_{k2^l} \bigr)^2
\bigr) \Biggr)^{1/2}
\\
&&\qquad = 2 \bigl(\E_0\bigl(S_{2^r}^2\bigr)
\bigr)^{1/2}+ 2 \sum_{l=0}^{r-1}
\Biggl( \sum_{k=1}^{2^{r-l}-1} Q^{k2^l}
\bigl( \bigl(\E_0(S_{2^l}) \bigr)^2 \bigr)
\Biggr)^{1/2} \qquad\mbox{$\P$-a.s.}\nonumber
\end{eqnarray}
\end{prop}

\begin{pf*}{Proof of Proposition \ref{propquenched2}}
Let $v\ge0$ be an integer, fixed for the moment. Let $r> v$.
Then we have
\[
\max_{1\le k \le2^r} |S_k|\le\max_{1\le s\le2^{r-v}}
|S_{s2^v }| + 2^v \max_{1\le j \le2^{r}}|X_j|.
\]
Let $K\ge1$, be fixed for the moment. We have
\[
\max_{1\le j \le2^r}|X_j|^2 \le K^2+ \sum_{j=1}^{2^{r}} |X_j|^2{\mathbf
1}_{\{|X_j|\ge K\}}.
\]

Hence, applying Proposition \ref{propMP} to the stationary
sequence $(S_{(k+1)2^v}-S_{k2^v})_{k\ge0}$
adapted to the filtration $(\F_{k2^v})_{k\ge0}$, we obtain (with the
convention that $S_0=0$)
\begin{eqnarray*}
\E_0 \Bigl(\max_{1\le i\le2^r } |S_i|^2
\Bigr) &\ll&4^v K^2 + 4^v \sum
_{j=1}^{2^r}Q^j \bigl(|X_0|^2{
\mathbf1}_{\{|X_0|\ge K\}
}\bigr)+ \E_0\bigl(S_{2^r}^2
\bigr)
\\
&&{}+ \Biggl(\sum_{l=0}^{r-v-1} \Biggl(\sum
_{k=1}^{2^{r-v-l}-1} Q^{k2^{l+v}} \bigl(\bigl(
\E_0(S_{2^{l+v}})\bigr)^2 \bigr)
\Biggr)^{1/2} \Biggr)^2
\\
&\ll&4^v K^2 + 4^v \sum
_{j=1}^{2^r}Q^j \bigl(|X_0|^2{
\mathbf1}_{\{|X_0|\ge K\}
}\bigr)+ \E_0\bigl(S_{2^r}^2
\bigr) \\
&&{}+ 2^{r} \Biggl(\sum_{l=0}^{r-v-1}
\frac{ (\M
_{l+v} ( (\E_0(S_{2^{l+v}}))^2 )
)^{1/2}}{2^{(l+v)/2}} \Biggr)^2.
\end{eqnarray*}
By assumption $\E_0(S_{2^r}^2)=o(2^r)$ $\P$-a.s.
By (\ref{ergthmQ}),
$(\sum_{j=1}^{2^r}Q^j (|X_0|^2{\mathbf1}_{\{|X_0|\ge K\}}))/\break 2^r \to
\E(|X_0|^2{\mathbf1}_{\{|X_0|\ge K\}} | {\mathcal I})$ $\P$-a.s.
Since $\|\E(|X_0|^2{\mathbf1}_{\{|X_0|\ge K\}} | {\mathcal I})\|_1
\leq\break  \| X_0^2{\mathbf1}_{\{|X_0|\ge K\}} \|_1 \rightarrow\infty$,
there exists a subsequence $(K_j)$ such that\break
$ \E(|X_0|^2{\mathbf1}_{\{|X_0|\ge K_j\}} | {\mathcal I}) \rightarrow0$
$\P$-a.s. as $j\to\infty$. Hence taking the $\limsup_r$
and letting $j\to\infty$, we obtain
\[
\limsup_r \frac{\E_0(\max_{1\le i\le2^r } |S_i|^2)}{2^r} \ll\biggl(\sum
_{l\ge v} \frac{ (\M_{l} ( (\E_0(S_{2^{l}}))^2 )
)^{1/2}}{2^{l/2}} \biggr)^2 \qquad\mbox{$
\P$-a.s.}
\]
To finish the proof, it suffices to prove that the random variable
defined by the series on the
right-hand side is $\P$-a.s. finite. But it is in $\LL^{2,w}$ since,
by~(\ref{MW}),
\begin{eqnarray*}
\biggl\|\sum_{l\ge0} \frac{ (\M_{l} ( (\E_0(S_{2^{l}}))^2 )
)^{1/2}}{2^{l/2}} \biggr\|_{2,w}
&\le&\sum_{l\ge0} \frac{\| ( \M_{l} ( (\E_0(S_{2^{l}}))^2 )
)^{1/2}\|_{2,w}}{2^{l/2}} \nonumber\\[-8pt]\\[-8pt]
&\ll&\sum
_{l\ge0} \frac{\|\E_0(S_{2^l})\|_2}{2^{l/2}} <\infty.
\end{eqnarray*}
\upqed
\end{pf*}

%
\begin{appendix}\label{app}
\section{\texorpdfstring{Proof of Theorem \lowercase{\protect\ref{theoergo}}}{Proof of Theorem 4.7}} \label{sectergodic}

We make the proof for $T$ Dunford--Schwartz and $f$
real-valued since the proof in the case where $f$ is $\B$-valued is identical,
replacing $|\cdot|$ with $| \cdot|_\B$ when necessary.

Write $U_n(f)=f+\cdots+ T^{n-1}f$. Since $\psi$ is monotonic, it
follows from the subadditivity of $(\|U_n(f)\|_p)$
(see, e.g., \cite{PU}, Lemma 2.7,
and \cite{MP}, equation (92)) that (\ref{MW1}) is equivalent to
\[
\sum_n \frac{\|f+\cdots+ T^{2^n-1}f\|_p} {\psi(2^n)
2^{n/p}}=\sum
_n \frac{\|U_{2^n}(f)\|_p} {\psi(2^n) 2^{n/p}} < \infty.
\]
We proceed now as in the proof of Proposition \ref{propquenched2};
namely, we consider dyadic blocs. Let us give the hints. Let $v\geq0$
be an integer. For $r>v$, write that
\[
\max_{1\le k \le2^r} \bigl|U_k(f)\bigr|\le\max_{1\le s\le2^{r-v}}
\bigl|U_{s2^v }(f)\bigr| + 2^v \max_{1\le j \le
2^{r}}\bigl|T^j
f\bigr|.
\]
Using Proposition \ref{prop} to take care of the first term in the
right-hand side, it follows that
\begin{eqnarray*}
\max_{1\le k \le2^r} \bigl|U_k(f)\bigr|&\le&2^v \max
_{1\le j \le2^{r}}\bigl|T^j f\bigr| \\
&&{}+ 2^{r/p}\sum
_{k\ge0} \frac{ [\M_{k +v}(|U_{2^{k
+v}}(f)|^p) ]^{1/p}}{2^{(k+v)/p}}.
\end{eqnarray*}
We finish the proof by using arguments developped in the proof of
Proposition~\ref{propquenched2}.

\section{Auxiliary results} \label{secttechres}
%
%
\begin{lem} \label{lmaabel}
Let $\B$ be a Banach space and $(a_n)_{n \geq1}$ a $\B$-valued
sequence. The following are equivalent:
\begin{longlist}[(ii)]
\item[(i)] the series $ \sum_{n \geq1} a_n $ converges;
\item[(ii)] $\lim_{n \rightarrow\infty}n \sum_{k \geq n}
(k+1)^{-1} a_k =0$ and the series $ \sum_{n \geq1} \sum_{k \geq n}
(k+1)^{-1} a_k $ converges.
\end{longlist}
\end{lem}
The proof is omitted since it follows from standard arguments based on
Abel summation by part.

The next lemma is Lemma 19 in Merlev\`ede, Peligrad and Peligrad \cite
{MPP}. In their paper, the lemma is stated with $\ell=0$ and with
${\mathcal H}=\R$, but with similar arguments as those in their proof,
it works for any nonnegative integer $\ell$ and for adapted
stationary sequences with values in a normed space by replacing the
absolute values by the corresponding norms.

%
%
\begin{lem}
\label{SUB}Let $p \geq1$ and let $(X_n)_{n\in\mathbb{Z}}$ be an
adapted stationary sequence in $\LL^p({\mathcal H})$ in the sense of
Notation \ref{notation1}. For every $\gamma>0$,
$n\geq1$ and any integer $\ell\geq0$,
\[
\frac{1}{n^{\gamma}}\max_{1\leq k\leq n}\bigl\|\E_{- \ell}%
(S_{k})\bigr\|_{p, {\mathcal H}}\leq2^{3\gamma+3} \sum
_{k=n+1}^{6n}\frac{1}{k^{\gamma+1}}\bigl\|
\E_{-\ell}(S_{k})\bigr\|_{p, {\mathcal H}}.
\]
\end{lem}

%
%
\begin{prop} \label{passMnSn}
Let $p \in[ 2,4]$ and let $(X_n)_{n\in\mathbb{Z}}$ be an adapted and
stationary sequence in $\LL^p$ in the sense of Notation \ref
{notation1}. Assume that
(\ref{MWpalpha}) holds. Then setting $M_n=\sum_{k=1}^{n}D\circ\theta
^k$ where $D$ is defined by (\ref{defdiffmart}),\vadjust{\goodbreak} the following
inequality holds: for any
nonnegative integers $r$ and $n$,
\begin{eqnarray*}
\bigl\|\E_{-r} \bigl(M_n^2\bigr) - \E
\bigl(M_n^2\bigr) \bigr\|_{p/2} &\ll&\bigl\|\E
_{-r}\bigl(S_n^2\bigr) - \E
\bigl(S_n^2\bigr) \bigr\|_{p/2} + \bigl\|
\E_{-r}\bigl(S^2_{2n}\bigr) - \E\bigl(
S^2_{2n} \bigr) \bigr\|_{p/2 }
\\
&&{} + n \biggl( \sum_{k \geq[n^{p/2}] } \frac{\|\E_0(S_k)\|
_{p}}{k^{1+1/p}}
\biggr)^2 + n \sum_{k \geq n }
\frac{\|\E
_{-n}(S_k)\|_{2}}{k^{3/2}}.
\end{eqnarray*}
In the statement of the proposition as well as in its proof,
the constants arising from the symbol $\ll$ are independent from $n$
and $r$.
\end{prop}
\begin{pf}
Setting $R_n = S_n - M_n$, we start with the following inequality:
%
%
\begin{eqnarray}
\label{deccarreprop} \bigl\|\E_{-r}\bigl(M_n^2
\bigr) - \E\bigl(M_n^2\bigr) \bigr\|_{p/2} &\ll&\bigl\|
\E_{-r}\bigl(S_n^2\bigr) - \E
\bigl(S_n^2\bigr) \bigr\|_{p/2} +2 \|
R_n \|^2_{p} \nonumber\\[-8pt]\\[-8pt]
&&{}+2 \bigl\|\E_{-r}
(S_nR_n) - \E(S_nR_n)
\bigr\|_{p/2}.\nonumber
\end{eqnarray}
Using Theorem \ref{theoMW} with $p \geq2$, we first get that
%
%
\begin{equation}
\label{b1Rn} \| R_n \|^2_{p} \ll n
\biggl( \sum_{k \geq[n^{p/2}] } \frac
{\|\E_0(S_k)\|_{p}}{k^{1+1/p}}
\biggr)^2.
\end{equation}
Now, starting from (\ref{dec}) and using the decompositions (\ref
{Vn2}), (\ref{Wn2}), (\ref{Vn}) and (\ref{calcul}) with $N=2n$, we
write that
%
%
\begin{equation}
\label{dec1Rndec} R_n = \E_0(S_n) +
\frac{\E_0(S_{2n} \circ\theta^n)}{2n+1} - \frac
{\E_n(S_{2n} \circ\theta^n)}{2n+1} -A_n -B_n,
\end{equation}
where
%
%
\begin{equation}
\label{decAn} A_n=2n \sum_{l\ge2n+1}
\frac{\E_{n}
(S_{l}\circ\theta^n) - \E_{0}
(S_{l}\circ\theta^n)}{l(l+1)}
\end{equation}
and
%
%
\begin{equation}
\label{decBn} B_n =\sum_{k\ge2n}\sum
_{l\ge k+1}\frac{\E_{n}(X_{l+n})-\E
_{0}(X_{l+n})}{l}.
\end{equation}
Notice first that
\begin{eqnarray*}
&& \biggl\| \E_{-r} \biggl( S_n \biggl( \E_0(S_n)
+\frac{\E_0(S_{2n} \circ\theta^n)}{2n+1} \biggr) \biggr) - \E\biggl(
S_n \biggl(
\E_0(S_n) +\frac{\E_0(S_{2n} \circ\theta
^n)}{2n+1} \biggr) \biggr)
\biggr\|_{p/2}
\\
&&\qquad\leq2 \biggl\| \E_{0} \biggl( S_n \biggl(
\E_0(S_n) +\frac{\E
_0(S_{2n} \circ\theta^n)}{2n+1} \biggr) \biggr)
\biggr\|_{p/2} \\
&&\qquad\leq2 \bigl\|\E_0 (S_n )
\bigr\|^2_{p} +2(2n+1)^{-1} \bigl\|\E_0
(S_n ) \bigr\|_{p} \bigl\|\E_0 (S_{2n}
) \bigr\|_{p},
\end{eqnarray*}
which combined with (\ref{applisub}) with $u_n=[n^{p/2}]$ implies that
%
%
\begin{eqnarray}
\label{b2Rn}
&&\biggl\| \E_{-r} \biggl( S_n \biggl(
\E_0(S_n) +\frac{\E_0(S_{2n}
\circ\theta^n)}{2n+1} \biggr) \biggr) \nonumber\\
&&\quad{}- \E
\biggl( S_n \biggl( \E_0(S_n) +
\frac{\E_0(S_{2n} \circ\theta
^n)}{2n+1} \biggr) \biggr) \biggr\|_{p/2}
\\
&&\qquad\ll n \biggl( \sum_{k\geq[n^{p/2}] } \frac{\|\E_0(S_k)\|
_{p}}{k^{1+1/p}}
\biggr)^2.\nonumber
\end{eqnarray}
Now writing that $S_{2n} \circ\theta^n = S_{2n} \circ\theta^n - S_n
\circ\theta^n + S_n \circ\theta^n$ and using the fact that $S_n$ is
${\mathcal F}_n$-measurable, we get
%
%
\begin{eqnarray}
\label{b3Rn*}
&& \biggl\| \E_{-r} \biggl( S_n \biggl(
\frac{\E_n(S_{2n} \circ\theta^n)}{2n+1} \biggr) \biggr) - \E\biggl(
S_n \biggl(
\frac{\E_n(S_{2n} \circ\theta^n)}{2n+1} \biggr) \biggr) \biggr\|_{p/2}
\nonumber\\
&&\qquad\leq n^{-1}\bigl\|\E_{-r} \bigl(S_n(S_{2n}-
S_n) \bigr) - \E\bigl( S_n (S_{2n}-
S_n )\bigr) \bigr\|_{p/2}\\
&&\qquad\quad{} +n^{-1} \bigl\|
\E_{-r} \bigl(S_n \E_n( S_{3n}
-S_{2n} ) \bigr)\bigr\|_{p/2}.\nonumber
\end{eqnarray}
Using the identity $2ab = (a+b)^2 - a^2 -b^2$ and the stationarity, we
first obtain that
%
%
\begin{eqnarray}
\label{b3Rn**}
&&2\bigl\|\E_{-r} \bigl(S_n(S_{2n}-
S_n) \bigr) - \E\bigl( S_n (S_{2n}-
S_n )\bigr) \bigr\|_{p/2}\nonumber\\
&&\qquad \leq 2\bigl\|\E_{-r}
\bigl(S^2_{n}\bigr) - \E\bigl( S^2_{n}
\bigr) \bigr\|_{p/2}
\\
&&\qquad\quad{} +\bigl\|\E_{-r} \bigl(S^2_{2n}\bigr) - \E\bigl(
S^2_{2n} \bigr) \bigr\|_{p/2}.\nonumber
\end{eqnarray}
To bound up the second term in (\ref{b3Rn*}), we write $C_n:=n^{-1} \E
_n( S_{3n} -S_{2n} )$, and we follow the lines of the proof of Theorem
2.3 in \cite{DDM}; see the display lines between their equations
(4.13) and (4.16). Hence we first write that
\begin{eqnarray*}
&&
\bigl\|\E_{-r} (S_nC_n ) \bigr\|_{p/2} \\
&&\qquad
\leq \bigl\|\E^{1/2}_{-r} \bigl(S_n^2
\bigr) \E^{1/2}_{-r}\bigl(C_n^2 \bigr)
\bigr\|_{p/2}
\\
&&\qquad\leq \bigl\|\bigl( \E_{-r} \bigl(S_n^2\bigr)-
\E\bigl(S_n^2\bigr) \bigr)^{1/2}
\E^{1/2}_{-r}\bigl(C_n^2 \bigr)
\bigr\|_{p/2}+ \bigl(\E\bigl(S_n^2\bigr)
\bigr)^{1/2}\bigl\|\E^{1/2}_{-r}\bigl(C_n^2
\bigr) \bigr\|_{p/2}
\\
&&\qquad \leq \bigl\|\E_{-r} \bigl(S_n^2\bigr)- \E
\bigl(S_n^2\bigr) \bigr\|_{p/2} + \bigl\|
C_n\bigr\|^2_{p} + \bigl(\E\bigl(S_n^2
\bigr) \bigr)^{1/2}\bigl\|\E^{1/2}_{-r}
\bigl(C_n^2 \bigr) \bigr\|_{p/2}.
\end{eqnarray*}
Notice that since (\ref{MWpalpha}) holds, by Theorem \ref{theoMW}, we
have in particular that
$
\| S_n \|_{2} = o(n^{1/2}) + \| M_n \|_{2}
$, implying that
%
%
\begin{equation}\label{majsnp}
\| S_n \|_{2} \ll n^{1/2}.
\end{equation}
Using (\ref{majsnp}) and the fact that the function $x \mapsto
|x|^{p/4}$ is concave, it follows that
%
%
\begin{equation}
\label{b44*Rn} \bigl\|\E_{-r} (S_nC_n )
\bigr\|_{p/2}\ll\bigl\|\E_{-r} \bigl(S_n^2
\bigr)- \E\bigl(S_n^2\bigr) \bigr\|_{p/2} + \|
C_n \|^2_{p} + n^{1/2} \|
C_n \|_{2}.
\end{equation}
By stationarity and using (\ref{applisub}) with $u_n=[n^{p/2}]$, we
get that
%
%
\begin{equation}
\label{b44*Rn*} \| C_n \|_{p}\ll n^{-1}
\bigl\|\E_{-n}( S_{n} ) \bigr\|_{p,{\mathcal H}} \ll
n^{-1/2}\sum_{k \geq[n^{p/2}] } \frac{\|\E
_0(S_k)\|_{p}}{k^{1+1/p}}.
\end{equation}
On the other hand, by using once again stationarity and Lemma \ref{SUB},
%
%
\begin{equation}
\label{b44*Rn**} \| C_n \|_{2}\ll n^{-1}
\bigl\|\E_{-n}( S_{n} ) \bigr\|_{2} \ll\sum
_{k \geq n } \frac{\|\E_{-n}(S_k)\|_{2}}{k^{2}}.
\end{equation}
Therefore starting from (\ref{b3Rn*}) and using (\ref{b3Rn**}),
(\ref{b44*Rn}), (\ref{b44*Rn*}) and (\ref{b44*Rn**}), we infer that
%
%
\begin{eqnarray}
\label{b3Rn}
&&\biggl\| \E_{-r} \biggl( S_n \biggl(
\frac{\E_n(S_{2n} \circ\theta_n)}{2n+1} \biggr) \biggr) - \E\biggl(
S_n \biggl(
\frac{\E_n(S_{2n} \circ\theta_n)}{2n+1} \biggr) \biggr) \biggr\|_{p/2} \nonumber\\
&&\qquad\ll
\bigl\|\E_{-r}
\bigl(S^2_{n}\bigr) - \E\bigl( S^2_{n}
\bigr) \bigr\|_{p/2}
\nonumber\\[-8pt]\\[-8pt]
&&\qquad\quad{} + n^{-1}\bigl\|\E_{-r}\bigl(S^2_{2n}
\bigr) - \E\bigl( S^2_{2n} \bigr) \bigr\|_{p/2} +
n^{-1} \biggl( \sum_{k \geq[n^{p/2}] }
\frac{\|\E_0(S_k)\|
_{p}}{k^{1+1/p}} \biggr)^2\nonumber\\
&&\qquad\quad{} + n^{1/2}\sum
_{k \geq n } \frac{\|\E
_{-n}(S_k)\|_{2}}{k^{2}}.\nonumber
\end{eqnarray}
We consider now the term $\|\E_{-r} (S_n A_n) - \E(S_n A_n) \|
_{p/2}$. With this aim, we first define
\[
\widetilde A_n= 2n \E_{n} \bigl(S_{n}\circ
\theta^n\bigr) \sum_{l\ge2n+1}
\frac{1}{ l(l+1)}.
\]
Since $S_n$ is ${\mathcal F}_n$-measurable,
\[
\bigl\|\E_{-r} (S_n \widetilde A_n) -
\E(S_n \widetilde A_n) \bigr\|_{p/2} \leq\bigl\|
\E_{-r} \bigl(S_n(S_{2n}- S_{n})
\bigr) - \E\bigl( S_n (S_{2n}- S_{n})\bigr)
\bigr\|_{p/2}.
\]
Using then the identity $2ab = (a+b)^2 - a^2 -b^2$ and stationarity, it
follows that
%
%
\begin{eqnarray}
\label{b24*Rn}&& 2 \bigl\|\E_{-r} (S_n \widetilde
A_n) - \E(S_n \widetilde A_n) \bigr\|
_{p/2} \nonumber\\[-8pt]\\[-8pt]
&&\qquad\leq2 \bigl\|\E_{-r}\bigl(S^2_n
\bigr) - \E\bigl(S_n^2 \bigr) \bigr\|_{p/2} + \bigl\|
\E_{-r}\bigl( S^2_{2n} \bigr)- \E\bigl(
S^2_{2n} \bigr)\bigr\|_{p/2}.\nonumber
\end{eqnarray}
Let now
\[
D_n:=n\sum_{k \geq2n+1} \frac{\E_{n}
(S_{k}\circ\theta^n) -\E_{n}
(S_{n}\circ\theta^n) }{k(k+1)}
\]
and notice that, by stationarity,
%
%
\begin{eqnarray}
\label{b34*Rn}
&&
\bigl\|\E_{-r} \bigl(S_n (A_n -
\widetilde A_n)\bigr) - \E\bigl(S_n (A_n-
\widetilde A_n)\bigr) \bigr\|_{p/2}
\nonumber\\[-8pt]\\[-8pt]
&&\qquad\ll n \bigl\|\E_0 (S_n) \bigr\|_p \sum
_{k \geq n+1} \frac{\|\E
_0 (S_k) \|_p}{k^2} + \bigl\|\E_{-r}
(S_n D_n )\bigr\|_{p/2}.\nonumber
\end{eqnarray}
Using (\ref{applisub}) with $u_n=n$, we first get that
\[
n \bigl\|\E_0 (S_n ) \bigr\|_{p} \sum
_{\ell\ge n+1} \frac{ \|\E
_0( S_{\ell} ) \|_{p} }{\ell^{2}} \ll n^2 \biggl( \sum
_{k \geq
n } \frac{\|\E_0(S_k)\|_{p}}{k^{2}} \biggr)^2.
\]
But, by using Lemma \ref{SUB} and the fact that $p \geq2$,
%
%
\begin{eqnarray}
\label{b3Rnprime}\quad n \sum_{k \geq n } \frac{\|\E_0(S_k)\|_{p}}{k^{2}}
& \leq& \max_{1
\leq k \leq[n^{p/2}]}\bigl\|\E_0(S_k)
\bigr\|_{p} + n \sum_{k \geq[n^{p/2}] } \frac{\|\E_0(S_k)\|_{p}}{k^{2}}
\nonumber\\[-8pt]\\[-8pt]
& \ll& n^{p/2} \sum_{k \geq[n^{p/2}] }
\frac{\|\E_0(S_k)\|
_{p}}{k^{2}} \ll n^{1/2} \sum_{k \geq[n^{p/2}] }
\frac{\|\E_0(S_k)\|
_{p}}{k^{1+1/p}}.\nonumber
\end{eqnarray}
Therefore,
%
%
\begin{equation}
\label{b334*Rn} n \bigl\|\E_0 (S_n ) \bigr\|_{p}
\sum_{\ell\ge n+1} \frac{ \|\E
_0( S_{\ell} ) \|_{p} }{\ell^{2}} \ll n \biggl( \sum
_{k \geq
[n^{p/2}] } \frac{\|\E_0(S_k)\|_{p}}{k^{1+1/p}} \biggr)^2.
\end{equation}
We bound now the second term in the right-hand side of (\ref{b34*Rn}).
Proceeding as to get (\ref{b44*Rn}), we infer that
%
%
\begin{equation}
\label{b44*Rnprime}\quad \bigl\|\E_{-r} (S_nD_n )
\bigr\|_{p/2}\ll\bigl\|\E_{-r} \bigl(S_n^2
\bigr)- \E\bigl(S_n^2\bigr) \bigr\|_{p/2} + \|
D_n \|^2_{p} + n^{1/2} \|
D_n \|_{2}.
\end{equation}
Stationarity and inequality (\ref{b3Rnprime}) imply that
%
%
\begin{equation}
\label{b54*Rn} \| D_n \|_{p} \ll n \sum
_{k \geq n} \frac{\|\E
_{-n}(S_k)\|_{p} }{k^2} \ll n^{1/2} \sum
_{k \geq[n^{p/2}] } \frac
{\|\E_0(S_k)\|_{p}}{k^{1+1/p}}.
\end{equation}
On the other hand, using once again stationarity,
%
%
\begin{equation}
\label{b64*Rn} \| D_n \|_{2} \ll n\sum
_{k \geq n} \frac{\|\E
_{-n}(S_k)\|_{2} }{k^2}.
\end{equation}
Overall, starting from (\ref{b34*Rn}) and considering the bounds
(\ref{b334*Rn}), (\ref{b44*Rnprime}), (\ref{b54*Rn}) and (\ref
{b64*Rn}), it follows that
%
%
\begin{eqnarray}
\label{b74*Rn}
&&
\bigl\|\E_{-r} \bigl(S_n (A_n
-\widetilde A_n)\bigr) - \E\bigl(S_n (A_n-
\widetilde A_n)\bigr) \bigr\|_{p/2} \nonumber\\
&&\qquad\ll\bigl\|\E_{-r}
\bigl(S^2_n \bigr) - \E\bigl(S_n^2
\bigr) \bigr\|_{p/2}
\\
&&\qquad\quad{} + n \biggl( \sum_{k \geq[n^{p/2}] } \frac{\|\E
_0(S_k)\|_{p}}{k^{1+1/p}}
\biggr)^2 + n^{3/2} \sum_{k \geq n}
\frac
{\|\E_{-n}(S_k)\|_{2} }{k^2}.\nonumber
\end{eqnarray}
We consider now the term $\|\E_{-r} (S_n B_n) - \E(S_n B_n) \|
_{p/2}$. Proceeding as to get (\ref{b44*Rn}), we infer that
%
%
\begin{eqnarray}
\label{b5Rn}
&&\bigl\|\E_{-r} (S_n B_n ) -
\E(S_nB_n ) \bigr\|_{p/2} \nonumber\\[-8pt]\\[-8pt]
&&\qquad\ll\bigl\|
\E_{-r} \bigl(S_n^2\bigr)- \E
\bigl(S_n^2\bigr) \bigr\|_{p/2} + \|
B_n \|^2_{p} + n^{1/2} \|
B_n \|_{2}.\nonumber
\end{eqnarray}
According to the bound (\ref{boundWn}) with $u_n = 2n$, followed by an
application of Lemma~\ref{lemmedalibor},
%
%
\begin{equation}
\label{b6Rn} \|B_n\|_{2}\ll n^{1/2} \sum
_{k\ge n}\sum_{m\ge1}
\frac{\|{\mathcal P}_{-k}(S_m)\|_{2}}{(m+k)^2} \ll n^{1/2} \sum_{k\ge n}
\frac{\|\E_{-n}(S_k)\|_{2}}{k^{3/2}}.
\end{equation}
To bound $\| B_n \|_{p} $, we use (\ref{dec1Rndec}). By
stationarity, we then infer that
\[
\| B_n \|_{p} \leq\| R_n
\|_{p} + 3 \bigl\|\E_0(S_n)
\bigr\|_{p} + 2 n \sum_{\ell\geq n+1}
\frac{\|\E_0(S_{\ell}) \|_{p}}{\ell^2}.
\]
Hence using Theorem \ref{theoMW} and inequality (\ref{applisub}) with
$u_n=n$, we get that
\[
\| B_n \|_{p} \ll n^{1/2}\sum
_{\ell\geq[n^{p/2}]} \frac{\|\E_0(S_{\ell}) \|_{p}}{\ell^{1+1/p}}
+ n \sum
_{\ell\geq n} \frac{\|\E_0(S_{\ell}) \|_{p}}{\ell^2},
\]
which together with (\ref{b3Rnprime}) implies that
%
%
\begin{equation}
\label{b7Rn} \| B_n \|_{p} \ll n^{1/2}\sum
_{\ell\geq[n^{p/2}]} \frac{\|\E_0(S_{\ell} )\|_{p}}{\ell^{1+1/p}}.
\end{equation}
Starting from (\ref{b5Rn}) and using (\ref{b6Rn}) and (\ref{b7Rn}),
we then obtain that
%
%
\begin{eqnarray}\label{b9Rn}
&&
\bigl\|\E_{-r} (S_n B_n ) -
\E(S_nB_n ) \bigr\|_{p/2} \nonumber\\
&&\qquad\ll \bigl\|
\E_{-r} \bigl(S_n^2\bigr)- \E
\bigl(S_n^2\bigr) \bigr\|_{p/2} + n \biggl( \sum
_{k \geq[n^{p/2}] } \frac{\|\E_0(S_k)\|
_{p}}{k^{1+1/p}} \biggr)^2
\\
&&\qquad\quad{} + n \sum_{k\ge n} \frac{\|\E_{-n}(S_k)\|_{2}}{k^{3/2}}.\nonumber
\end{eqnarray}
Taking into account the decomposition (\ref{dec1Rndec}) together with
the bounds (\ref{b2Rn}), (\ref{b3Rn}), (\ref{b24*Rn}), (\ref
{b74*Rn}) and (\ref{b9Rn}), we then derive that
%
%
\begin{eqnarray}
\label{b10Rn}
&&\bigl\|\E_{-r} (S_n R_n ) -
\E(S_nR_n ) \bigr\|_{p/2} \nonumber\\
&&\qquad\ll\bigl\|
\E_{-r} \bigl(S_n^2\bigr) - \E
\bigl(S_n^2\bigr) \bigr\|_{p/2} +\bigl\|
\E_{-r}\bigl(S^2_{2n}\bigr) - \E\bigl(
S^2_{2n} \bigr) \bigr\|_{p/2}
\\
&&\qquad\quad{} + n \biggl( \sum_{k \geq[n^{p/2}] } \frac{\|\E_0(S_k)\|
_{p}}{k^{1+1/p}}
\biggr)^2 + n \sum_{k \geq n }
\frac{\|\E
_{-n}(S_k)\|_{2}}{k^{3/2}}.\nonumber
\end{eqnarray}
Starting from (\ref{deccarreprop}) and considering the inequalities
(\ref{b1Rn}) and (\ref{b10Rn}), the proposition follows.
\end{pf}
\end{appendix}

\section*{Acknowledgments}

C. Cuny would like to thank Dalibor Voln\'{y} for helpful discussions
while he was visiting the university of Rouen. The authors are also
indebted to the referee for carefully reading the manuscript.


%
%

\printaddresses

\end{document}